\documentclass{article}
\usepackage{amsfonts}
\usepackage{amsmath,amsthm}

\usepackage[all]{xy}
\usepackage{graphicx}

\usepackage{amssymb}
\usepackage{latexsym}

\DeclareMathAlphabet\EuFrak{U}{euf}{m}{n}	
\SetMathAlphabet\EuFrak{bold}{U}{euf}{b}{n}	

\newcommand{\ra}{\rightarrow}

\newcommand{\hra}{\hookrightarrow}

\newcommand{\wa}{\widehat}

\newcommand{\sC}{{\it C*}-}
\newcommand{\bC}{{\mathbb C}}

\newcommand{\bT}{{\mathbb T}}

\newcommand{\bZ}{{\mathbb Z}}
\newcommand{\bM}{{\mathbb M}}
\newcommand{\bN}{{\mathbb N}}
\newcommand{\bP}{{\mathbb P}}
\newcommand{\bS}{{\mathbb S}}

\newcommand{\ud}{{{\mathbb U}(d)}}
\newcommand{\sud}{{{\mathbb {SU}}(d)}}

\newcommand{\mA}{\mathcal A}
\newcommand{\mB}{\mathcal B}
\newcommand{\mcB}{\mathcal B}

\newcommand{\mE}{\mathcal E}
\newcommand{\mF}{\mathcal F}

\newcommand{\mL}{\mathcal L}
\newcommand{\mM}{\mathcal M}
\newcommand{\mN}{\mathcal N}
\newcommand{\mO}{\mathcal O}

\newcommand{\mR}{\mathcal R}
\newcommand{\mT}{\mathcal T}

\newcommand{\mZ}{\mathcal Z}

\newcommand{\mcG}{\mathcal G}
\newcommand{\mG}{G}

\newcommand{\vL}{L}
\newcommand{\vM}{M}

\newcommand{\tend}{{\bf end}\mA}

\newcommand{\spzro}{ X^\rho }

\newcommand{\ii}{\iota,\iota}
\newcommand{\rr}{\rho,\rho}

\newcommand{\mrs}{\mM^r,\mM^s}
\newcommand{\ers}{\mE^r,\mE^s}
\newcommand{\hrs}{H^r , H^s}
\newcommand{\rhors}{\rho^r , \rho^s}

\newcommand{\wE}{\wa {\mE}}
\newcommand{\wL}{\wa {\mL}}

\newcommand{\mcSUE}{\mathcal {SU} \mE}
\newcommand{\mcUE}{\mathcal {U} \mE}

\newcommand{\mSUE}{ {{\bf {SU}} \mE} }
\newcommand{\mUE}{ { {\bf U} \mE} }

\newcommand{\mSG}{ S \mG } 
\newcommand{\mUM}{{{\bf U} \mM} }

\newcommand{\cpn}{ \mA \rtimes_\rho \bN }
\newcommand{\cpen}{ \mA \rtimes_\rho^\mE \bN }

\newcommand{\coe}{\mO_{\mE}}

\newcommand{\com}{\mO_{\mM}}
\newcommand{\toe}{\mT_{\mE}}
\newcommand{\bm}{\mcB_{\mM}}

\newcommand{\cog}{\mO_{\mG}}
\newcommand{\cmg}{\mO^{\mG}_{\mM}}

\newcommand{\aog}{{^0\mathcal O_{\mG}}}

\newtheorem{thm}{Theorem}[section]
\newtheorem{cor}[thm]{Corollary}
\newtheorem{lem}[thm]{Lemma}
\newtheorem{prop}[thm]{Proposition}

\newtheorem{defn}[thm]{Definition}

\theoremstyle{definition}
\newtheorem{ex}{Example}[section]

\theoremstyle{remark}
\newtheorem{rem}{Remark}[section]

\numberwithin{equation}{section}

\begin{document}

\author{{\sf Ezio Vasselli}
                         \\{\it Dipartimento di Matematica}
                         \\{\it University of Rome "La Sapienza"}
			 \\{\it P.le Aldo Moro, 2 - 00185 Roma - Italy }
			 \\{\it c/o Sergio Doplicher }
                         \\{\sf vasselli@mat.uniroma2.it}}

\title{ Crossed products by endomorphisms,\\vector bundles\\and\\group duality }
\maketitle

\begin{abstract}
We construct the crossed product $\cpen$ of a $C(X)$-algebra $\mA$ by an endomorphism $\rho$, in such a way that $\rho$ becomes induced by the bimodule $\wE$ of continuous sections of a vector bundle $\mE \ra X$. Some motivating examples for such a construction are given. Furthermore, we study the \sC algebra of $\mG$-invariant elements of the Cuntz-Pimsner algebra $\coe$ associated with $\wE$, where $\mG$ is a (noncompact, in general) group acting on $\mE$. In particular, the \sC algebra of invariant elements w.r.t. the action of the group of special unitaries of $\mE$ is a crossed product in the above sense. We also study the analogous construction on certain Hilbert bimodules, called 'noncommutative pullbacks'.

\bigskip

\noindent {\em AMS Subj. Class.:} 46L05, 46L08, 22D35.

\noindent {\em Keywords:} Cuntz-Pimsner algebra; crossed product; continuous bundle; vector bundle; Hilbert bimodule.

\end{abstract}


\section{Introduction.}
\label{intro}

Let $d \in \bN$, $\ud$ the unitary group of $d \times d$-matrices, $\mO_d$ the Cuntz algebra generated by a set $\left\{ \psi_h \right\}_{h=1}^d$ of isometries with relations $\psi_h^* \psi_k = \delta_{hk}1$, $\sum_k \psi_k \psi_k^* = 1$. It is well-known that for every closed group $G \subseteq \ud$ there is a $G$-action by automorphisms on $\mO_d$, obtained in the following way: let $H \subset \mO_d$ be the vector space spanned by $\left\{ \psi_h \right\}$; then $H$, endowed with the scalar product $\psi , \psi' \mapsto \psi^* \psi' \in \bC 1$, $\psi , \psi' \in H$, is isomorphic to the standard rank $d$ Hilbert space $\bC^d$, so that there is an action $G \times H \ra H$ by unitary operators, corresponding to the defining representation of $G$ over $\bC^d$. By universality of the Cuntz algebra, such an action extends to an automorphic action 
\[
G \ra {\bf aut} \mO_d \ \ , \ \ 
G \ni g \mapsto \wa g \in {\bf aut} \mO_d \ .
\]
\noindent We denote by $\mO_G \subseteq \mO_d$ the fixed-point algebra. Important cases are the \sC algebras $\mO_\ud$ and $\mO_\sud$, where $\sud$ is the special unitary group: in particular, it turns out that $\mO_\sud$ is the crossed product of $\mO_\ud$ by a certain shift endomorphism (\cite[Lemma 3.8]{DR87}).

Now, $\mO_d$ is equipped with the {\em canonical endomorphism}
\begin{equation}
\label{def_eod}
\sigma \in {\bf end}\mO_d \ \ : \ \  
\sigma (t) := \sum_k \psi_k t \psi_k^* \ \ , \ \ t \in \mO_d \ ;
\end{equation}
\noindent since $\sigma \circ \wa g = \wa g \circ \sigma$, $g \in G$, we find that $\sigma$ restricts to an endomorphism $\sigma_G \in {\bf end} \mO_G$. The pair $(\mO_G , \sigma_G)$ can be used to prove a version of the Tannaka duality; in order to expose this fact, we recall the following general construction.

Let $\mA$ be a \sC algebra, $\tend$ the set of endomorphisms of $\mA$, $\rho \in \tend$. It is possible to associate with $\rho$ a tensor category, say $\wa \rho$, with objects the endomorphisms $\rho^r$, $r \in \bN$ (for $r = 0$, we define $\rho^0$ as the identity automorphism $\iota$); the spaces of arrows and the tensor structure are given by
\[
\left\{ 
\begin{array}{l}
(\rho^r , \rho^s) := \left\{ t \in \mA : t \rho^r (a) = \rho^s (a)t  \ , \ a \in \mA \right\} \ , \\

\rho^r , \rho^s \mapsto \rho^r \times \rho^s := \rho^{r+s}  \ , \\ 

t,t' \mapsto t \times t' := t \rho^r (t') = \rho^s (t')t 
\in (\rho^{r+r'} , \rho^{s+s'}) \ ,
\end{array}
\right.
\]
\noindent $t \in (\rho^r,\rho^s)$, $t' \in (\rho^{r'} ,\rho^{s'})$. Note that $(\ii) = \mA \cap \mA'$.

We now return to the pair $(\mO_G , \sigma_G)$. Let $\wa G$ denote the category of tensor powers of the defining representation of $G$. Then, it is proved that there is an isomorphism of tensor categories $\wa \sigma_G \simeq \wa G$ (\cite[Thm.3.5]{DR87}); moreover, $G$ can be characterized as the stabilizer of $\mO_G$ in $\mO_d$ (\cite[Cor.3.3]{DR87}). Thus, we are able to reconstruct the group $G$ and the category $\wa G$ by analyzing the algebraic properties of the \sC dynamical system $(\mO_G , \sigma_G)$, together with the inclusion $\mO_G \subseteq \mO_d$.

\bigskip

The above construction is an important step towards the Doplicher-Roberts duality for compact groups (\cite{DR89}), characterizing certain abstract tensor categories as group duals. In the present work, we start the program of generalizing such a duality to the case in which the set of arrows of the identity object does not reduce to the complex numbers. As a first step, we generalize the above 'Tannaka duality' for closed subgroups of $\ud$ to the case of certain noncompact groups naturally acting on vector bundles. Instead of the Cuntz algebra, we will make use of the Cuntz-Pimsner algebra associated with the module of continuous sections of a vector bundle $\mE \ra X$. The tensor categories that we consider are characterized by the property $(\ii) = C(X)$, where $\iota$ is the identity object. In a future paper, we will proceed with our program by constructing a crossed product of a \sC algebra by an endomorphism satisfying a weaker version of permutation symmetry and special conjugate property w.r.t. the analogous notions due to Doplicher and Roberts (\cite{DR89A}). This will generalize \cite[Thm.4.1]{DR89A} (recovered as the trivial-centre case), and \cite[Thm.2.4]{BL01} for the case of a tensor category generated by a single object (in our case, the generating "admissible" DR-category in the sense of the above-cited reference is replaced by a more general category with $(\ii) \neq \bC$). The groups we recover are in general non-compact.

The present paper is organized as follows:

\bigskip

In Section \ref{cp_vec_bundle} we construct the crossed product of a \sC algebra $\mA$ by an endomorphism $\rho$, in the following way: we consider the centre $ZM(\mA)$ of the multiplier algebra, and the unital abelian \sC algebra 
\begin{equation}
\label{def_zro}
C(\spzro):= \left\{ 
f \in ZM(\mA) \\ : \ \rho (fa) = f \rho(a) \ , \ a \in \mA  
\right\} \ \ .
\end{equation}
\noindent Then, we consider a vector bundle $\mE \ra \spzro$, and construct a crossed product in such a way that $\rho$ becomes induced (in the sense of Def.\ref{def_inn}) by the module of continuous sections of $\mE$ (Prop.\ref{prop_cp}). Our construction, made in the context of $C(X)$-algebras (\cite{Kas88,Bla96}), includes as particular cases (i.e. trivial vector bundles) the crossed products by Stacey \cite{Sta93}, Cuntz \cite{Cun77} and Paschke \cite{Pas80}. Motivating examples arise in the context of continuous trace \sC algebras (Ex.\ref{ex_ct}), and Cuntz-Pimsner algebras associated with vector bundles (Ex.\ref{ex_cp}).

\bigskip

In Section \ref{cavb_gb} we introduce the \sC algebra $\cog$ associated with a closed group $\mG$ of unitaries acting on a vector bundle $\mE \ra X$. We associate with $\mG$ a bundle $\mcG \ra X$, called the {\em spectral bundle} of $\mG$ (Def.\ref{def_group_bundle}): roughly speaking, $\mG$ is characterized as a 'total subset' of continuous sections of $\mcG$. We introduce the notion of {\em dual} of $\mG$, as the analogue of the category of tensor powers of the defining representation of a compact Lie group (Def.\ref{def_g_dual}); then, we study the analogue of the above mentioned construction for $(\cog , \sigma_\mG)$ (Cor.\ref{cor_itc}, Prop.\ref{dual_group_bundle}). An interesting fact is that non-isomorphic groups may have isomorphic duals; there are two orders of reason for this phenomenon:

\begin{itemize}

\item at a first level, it is verified that different groups have the same dual if and only if the associated spectral bundles coincide (Lemma \ref{lem_gsg});

\item at a deeper level, a dual may have several embeddings into the tensor category of vector bundles; different embeddings give rise to non-isomorphic groups (Cor.\ref{no_uni}, Ex.\ref{ex_sue}).

\end{itemize}

\noindent  In the important case in which $\mG$ is the group of special unitaries of $\mE$, then $\cog$ is a crossed product in the sense of the previous section (Prop.\ref{lem22}): this generalizes the analogous result for the \sC algebra $\mO_\sud$.

\bigskip

In order for future applications, in Section \ref{bim_pul} we introduce and describe a class of Hilbert bimodules corresponding to noncommutative pullbacks of vector bundles, in the sense of Def.\ref{def_bp}. The interest in such bimodules arises from the fact that the associated Cuntz-Pimsner algebras exhibit canonical endomorphisms with weaker properties w.r.t. the case of vector bundles (in particular, they do not have permutation symmetry in the sense of \cite[\S 4]{DR89A}). Vector bundle deformations in the sense of \cite{Wal02}, and vector $\mA$-bundles in the sense of Mishchenko (\cite{MF80}) supply examples of noncommutative pullbacks. We study group actions over noncommutative pullbacks and the associated Cuntz-Pimsner algebras (Prop.\ref{cross_bp}, Cor.\ref{cross_bp1}), generalizing some results of the previous section.

\section{Keywords.}
\label{preli}

The main references for {\em tensor} \sC {\em categories} are \cite{DR89,LR97}, to which we refere for the notion of {\em symmetry}. Tensor \sC categories of \sC algebra endomorphisms, and their relationship with the duality theory, are studied in \cite{DR87} in the setting of the Cuntz algebra, and more in general in \cite{DR89A}. We refere to these papers for the notions of {\em permutation symmetry} and {\em special conjugate property}. If $G \subseteq \ud$, then $\sigma_G \in {\bf end} \cog$ has permutation symmetry (\cite[\S 4]{DR87}), if $G \subseteq \sud$, then $\sigma_G$ satisfies also the special conjugate property (\cite[Lemma 2.2]{DR87}).

\bigskip

Let $X$ be a locally compact Hausdorff space. A {\em continuous bundle} of \sC algebras over $X$ is a \sC algebra $\mA$, equipped with a faithful family of epimorphisms $\left\{ \pi_x : \mA \ra \mA_x \right\}_{x \in X}$ such that for every $a \in \mA$ the norm function $\left\{ x \mapsto \left\| \pi_x (a) \right\| \right\}$ belongs to $C_0(X)$, and $\mA$ is a nondegenerate $C_0(X)$-bimodule w.r.t. pointwise multiplication $f , a \mapsto \left\{ f(x) \cdot \pi_x (a) \right\}$, $f \in C_0(X)$. The term \sC {\em algebra bundle} will be also used. Standard references for \sC algebra bundles (and the related notion of {\em continuous field}) are \cite[\S 10]{Dix},\cite{DD63,KW95}, to which we refere for the notions of {\em restriction} and {\em local triviality}. A $C_0(X)$-{\em algebra} is a \sC algebra $\mA$, equipped with a nondegenerate morphism from $C_0(X)$ into the centre of the multiplier algebra $M(\mA)$; in the sequel, we will identify elements of $C_0(X)$ with their image in $M(\mA)$. $C_0(X)$-algebra morphisms are by definition \sC algebra morphisms which are equivariant w.r.t. the $C_0(X)$-module actions. $C_0(X)$-algebras were introduced by Kasparov in \cite[Def.1.5]{Kas88}, and characterized as 'upper semicontinuous bundles' over $X$ in \cite{Bla96,Nil96}. As can be expected, continuous bundles of \sC algebras are $C_0(X)$-algebras. Every $C_0(X)$-algebra $\mA$ can be faithfully represented over a suitable Hilbert $C_0(X)''$-module, where $C_0(X)''$ denotes the enveloping algebra (see \cite[\S 1.6]{Kas88}).

\bigskip

Let $\mA$ be a \sC algebra, $M(\mA)$ the multiplier algebra. The Cuntz-Pimsner algebra ({\em CP-algebra}, in the sequel) associated with a Hilbert $\mA$-bimodule $\mM$ has been introduced in \cite{Pim97}. In \cite[\S 3]{DPZ97}, the CP-algebras are defined by using a universal construction, essentially the one exposed in \cite[\S 4]{DR89},\cite[\S 5]{LR97}. The construction is (briefly) the following: we consider the Banach $\mA$-bimodules $(\mrs)$ of right $\mA$-module operators from the (internal) tensor power $\mM^r$ into $\mM^s$, $r,s \in \bN$; in particular, we define $\mM^0 := \mA$, so that $(\mA , \mA) \simeq M(\mA)$. Note that $(\mM , \mM)$ is the \sC algebra of bounded right $\mA$-module operators of $\mM$. If $1 \in (\mM , \mM)$ is the identity operator, we embed $(\mrs)$ into $(\mM^{r+1}, \mM^{s+1})$ by tensoring on the right by $1$. Then, we consider the inductive limit $\com^k := \cdots (\mrs) \hra (\mM^{r+1}, \mM^{s+1}) \hra \cdots$, where $k := s - r \in \bZ$. If $t \in \com^h$, $t' \in \com^k$, we can define a product by assuming $t \in (\mM^{r}, \mM^{r+h})$, $t' \in (\mM^{r+h}, \mM^{r+h+k})$, and by composing $t' \cdot t$. In the same way, the operation of assigning the adjoint operator induces an involution $* : \com^k \ra \com^{-k}$. We define $\ ^0\com := \sum^\oplus_k \com^k$; the above considerations imply that $\ ^0\com$ is a $*$-algebra. It can be proved that there exist a unique \sC norm on $\ ^0\com$ such that the {\em circle action} $z , t \mapsto z^k t$, $z \in \bT$, $t \in \com^k$ extends to an isometric action. The closure w.r.t. such \sC norm is a \sC algebra $\com$, coinciding with the CP-algebra if $\mA$ is unital and $\mM$ is finitely generated (see also \cite{KPW98,KPW01} for a detailed study about the structure of CP-algebras). If $\mA = \bC$ and $\mM$ has finite rank $d \in \bN$, then $\com$ is the Cuntz algebra $\mO_d$.

\bigskip

Let $X$ be a compact Hausdorff space, $\mE \ra X$ a vector bundle. We denote by $\wE$ the finitely generated Hilbert $C(X)$-bimodule of continuous sections of $\mE$, endowed with coinciding left and right $C(X)$-module actions, and by $\coe$ the corresponding CP-algebra. Basic properties of $\coe$ are established in \cite[\S 4]{Vas}; we briefly recall them for the reader's convenience. According to the above mentioned procedure, $\coe$ is generated by the Banach $C(X)$-bimodules $(\wE^r , \wE^s)$, $r,s \in \bN$, which are identified (by the Serre-Swan theorem) with the the sets $(\ers)$ of vector bundle morphisms from the tensor power $\mE^r$ into $\mE^s$. In particular, note that $(\mE , \mE) \simeq ( \wE  , \wE )$ is the \sC algebra of bounded $C(X)$-module operators of $\wE$, and that there is a natural identification $\wE \simeq (\iota,\mE)$, where $\iota := \mE^0 := X \times \bC$. If $\left\{ \psi_l \right\}_l$ is any finite set of generators for $\wE$, then $\coe$ can be described in terms of generators and relations (as by \cite[\S 3]{Pim97}):
\begin{equation}
\label{genrel}
f  \psi_l = \psi_l  f \ \ , \ \
{\psi_l^*}  \psi_m = \left\langle \psi_l , \psi_m \right\rangle \ \ , \ \
\sum \limits_l \psi_l \psi_l^* = 1  \ \ ;
\end{equation}
\noindent here $f \in C(X)$ and $\left\langle \cdot , \cdot \right\rangle$ denotes the $C(X)$-valued scalar product on $\wE$. By removing the third of (\ref{genrel}), we get the Toeplitz algebra $\toe$ (in the sense of \cite{Pim97}). $\coe$ is a locally trivial continuous bundle of Cuntz algebras over $X$. This fact is proved in \cite[\S 4]{Vas}, but it was previously discovered by Roberts (\cite{Rob}). $\coe$ is endowed with a canonical shift endomorphism $\sigma$, coinciding with (\ref{def_eod}) in the case of the Cuntz algebra, in the following way: if $t \in (\ers) \subset \coe$, then  
\begin{equation}
\label{def_sigma}
\sigma (t) := 1 \otimes t \ \in (\mE^{r+1},\mE^{s+1}) \ \ ,
\end{equation}
\noindent where $1$ is the identity on $\mE$. Note that $\sigma$ is a $C(X)$-module map (i.e., $\sigma$ is the identity on $C(X) \subset \coe$). In the present paper, we will also adopt the notation $\psi_\vL := \psi_{l_1} \cdots \psi_{l_r}  \in  (\iota , \mE^r)$, $\vL = ( l_1 , \ldots , l_r )$, so that an element of $(\ers) \subset \coe$ is a linear combination with coefficients in $C(X)$ of terms of the type $\psi_\vM \psi_\vL^*$, $\vM := ( m_1 , \ldots , m_s )$; note that $\left\{ \psi_\vL \right\}$ is a set of generators for $\wa {\mE^r} \simeq (\iota , \mE^r)$. If $\psi , \psi' \in (\iota , \mE^r)$, $r \in \bN$, then the relation $\psi^* \psi' = \left \langle \psi ,\psi' \right \rangle \in C(X)$ holds in the CP-algebra, where $\left \langle \cdot ,\cdot \right \rangle$ is the $C(X)$-valued scalar product defined on the module of continuous sections of $\mE^r$.

\section{Crossed Products by Endomorphisms and Vector Bundles.}
\label{cp_vec_bundle}

We start the present section by exposing some elementary properties of Hilbert bimodules in \sC algebras (see \cite{DPZ97} for details). Let $\mB \subseteq \mA$ be an inclusion of \sC algebras; a closed vector space $\mM \subset \mA$ is said {\em Hilbert} $\mB$-{\em bimodule in} $\mA$ if it is stable for left and right multiplication by elements of $\mB$, and if $\psi^* \psi' \in \mB$, $\psi , \psi' \in \mM$; in such a way, we have a map
\[
\mM \times \mM \ra \mB \ \ ,\ \ \psi , \psi' \mapsto \psi^* \psi' 
\]
\noindent defining a $\mB$-valued scalar product on $\mM$. 

\bigskip

$\mM$ is {\em finitely generated} if there is a finite set $\left\{ \psi_l \right\} \subset \mM$ such that $\sum_l \psi_l \psi_l^* \psi = \psi$ for every $\psi \in \mM$ (note that $\psi_l^* \psi \in \mB$, thus $\psi$ is a linear combination with coefficients in $\mB$ of the generators). It is trivial to verify that $P := \sum_l \psi_l \psi_l^*$ is a projection, and that it does not depend on the choice of the generators. We call $P$ the {\em support} of $\mM$. Note that for every $r \in \bN$, there is a natural identification of $\mM^r :=$ ${\mathrm{span}}\left\{ \psi_1 \cdots \psi_r , \psi_i \in \mM , i = 1, \ldots , r \right\}$ with the $r$-fold internal tensor power $\mM \otimes_\mB \cdots \otimes_\mB \mM$. With the same argument, we have an identification $(\mrs) = \mM^s (\mM^*)^r :=$ ${\mathrm{span}} \left\{ \psi' \psi^* , \psi' \in \mM^s ,  \psi \in \mM^r  \right\}$. If $\mB \subseteq \mA \cap \mA'$, then an endomorphism $\sigma_\mM$ is defined on $\mA$: 
\begin{equation}
\label{12}
\sigma_\mM (a) := \sum_l \psi_l a \psi^*_l \ \ , \ \ a \in \mA \ .
\end{equation}
\noindent $\sigma_\mM$ does not depend on the choice of the generators, thus we say that $\sigma_\mM$ is the {\em inner endomorphism induced by} $\mM$. If $\mA$ has identity $1$, then $P = \sigma_\mM (1)$. Note that if $\mB = \mA \cap \mA'$, then 
\[
\mM = 
(\iota , \sigma_\mM) := 
\left\{ \psi \in \mA : \psi a = \sigma_\mM (a) \psi \ , \  a \in \mA \right\}
\ .
\]
\noindent In fact, if $\psi \in \mM$ then $\sigma_\mM (a) \psi = \sum_l \psi_l a (\psi_l^* \psi) = P \psi a = \psi a$, so that $\mM \subseteq (\iota , \sigma_\mM)$; viceversa, if $\psi' \in (\iota , \sigma_\mM)$ then $\psi^* \psi' \in (\ii) = \mB$ for every $\psi \in \mM$, thus $\psi' = \psi' \cdot 1 = \sigma_\mM (1) \psi' = \sum_l \psi_l (\psi_l^* \psi') \in \mM$.

\begin{defn}
\label{def_inn}
Let $\mA$ be a \sC algebra with centre $\mZ$, $\rho$ an endomorphism of $\mA$. $\rho$ is said {\bf inner} if there is a finitely generated Hilbert $\mZ$-bimodule $\mM \subset \mA$ such that $\rho = \sigma_\mM$.
\end{defn}

\begin{ex} {\em
Let $\mA$ be a \sC algebra with identity $1$ and centre $\mZ$, $v \in \mA$ a partial isometry with $v^* v = 1$, $v v^* = P$. Then, $\mM := \mZ \cdot v$ is a Hilbert $\mZ$-bimodule in $\mA$ with support $P$, (finitely) generated by $v$. An inner endomorphism $\sigma_\mM \in \tend$, $\sigma_\mM (a) := v a v^*$, $a \in \mA$, is induced.
} \end{ex}

\begin{ex} {\em
\label{ex_ivb}
Let $\mE \ra X$ be a vector bundle over a compact Hausdorff space; then, the module $\wE$ of continuous sections of $\mE$ is finitely generated by a set $\left\{ \psi_l \right\}$. We consider the canonical endomorphism $\sigma$ defined on the CP-algebra $\coe$ (see \S \ref{preli}). By \cite[Prop.4.2]{Vas} we find that $\sigma$ is inner, with $(\iota,\sigma) = (\iota , \mE) \simeq \wE$; so that, $\sigma (t) = \sum_l \psi_l t \psi_l^*$, $t \in \coe$. This also implies $(\ers) = (\sigma^r , \sigma^s)$, $r,s \in \bN$. In fact, $(\ers)$ is spanned as a vector space by the elementary tensors $\psi_\vL \psi_\vM^*$, that belong to $(\sigma^r , \sigma^s)$, so that $(\ers) \subseteq (\sigma^r , \sigma^s)$; viceversa, if $t \in (\sigma^r , \sigma^s)$ then $t_{LM} := \psi_\vL^* t \psi_\vM \in (\ii) = C(X)$, and $t = \sum_{LM} \psi_\vL t_{LM} \psi_\vM^* \in (\ers)$. In particular, we find that $(\iota , \sigma^s ) \simeq (\iota , \mE^s) \simeq \wa {\mE^s}$, so that $\sigma^s$ is inner.
} \end{ex}

\bigskip

Let $\rho$ be an endomorphism of a \sC algebra $\mA$. Then $\rho$ induces a natural structure of $C(X)$-algebra on $\mA$. In fact, if $M(\mA)$ is the multiplier algebra with centre $ZM(\mA)$, we consider the unital, abelian \sC algebra $C(\spzro) \subseteq ZM(\mA)$ defined in (\ref{def_zro}); since the identity $1_{M(\mA)} \in M(\mA)$ belongs to $C(\spzro)$, we conclude that $\mA$ is a $C(\spzro)$-algebra. Furthermore, by definition $\rho$ is a $C(\spzro)$-endomorphism.

\bigskip

Let now $X$ be a compact Hausdorff space, $\mA$ a $C(X)$-algebra, $\rho$ a $C(X)$-endomorphism of $\mA$. We consider $d \in \bN$ and a rank $d$ vector bundle $\mE \ra X$, so that $\wE$ is finitely generated as a Hilbert $C(X)$-module. We want to construct a crossed product of $\mA$ by $\rho$, in such a way that $\rho$ becomes induced by $\wE$ as by (\ref{12}).

\begin{defn}
A covariant representation of $(\mA,\rho)$ with rank $\mE$ is a pair $(\pi , \wE_\pi)$, where 

\begin{itemize}

\item $\pi : M(\mA) \ra L(\mM_\pi)$ is a unital $C(X)$-representation of the multiplier algebra $M(\mA)$ over a Hilbert $C(X)''$-module $\mM_\pi$;

\item $\wE_\pi$ is a Hilbert $C(X)$-bimodule in $L(\mM_\pi)$, isomorphic to $\wE$; 

\item $\pi \circ \rho = \sigma_{\pi} \circ \pi$, where $\sigma_\pi$ is the inner endomorphism induced by $\wE_\pi$ on $L (\mM_\pi)$ as by (\ref{12}).

\end{itemize}
\end{defn}

If $\mE$ is trivial (i.e. $\mE = X \times \bC^d$), we can find a set $\left\{ \psi_h \right\}_{h=1}^d$ of orthonormal sections generating $\wE$, i.e. $\left \langle \psi_h , \psi_k \right \rangle = \delta_{hk} 1$, where $\delta_{hk}$ is the Kroneker symbol. By the isomorphism $\wE_\pi \simeq \wE$, we obtain that $\wE_\pi$ is generated by a set $\left\{ \psi_{h,\pi} \right\}_{h=1}^d$ of orthogonal isometries such that $\psi_{h,\pi}^* \psi_{k,\pi} = \delta_{hk} 1$. Thus, by (\ref{12}) we find that $\sigma_\pi$ is induced by $\left\{ \psi_{h,\pi} \right\}_{h=1}^d$, as for the notion of covariant representation considered by Stacey (\cite[\S 2]{Sta93}). For example, if $\mE \simeq X \times \bC$, we find that $\sigma_\pi (\cdot) = v \cdot v^*$ for some partial isometry $v \in L(\mM_\pi)$.

We now prove the existence of covariant representations. For this purpose, let us denote by $\coe$ the CP-algebra of $\wE$, and by $\toe$ the corresponding Toeplitz algebra. By universality of the Toeplitz-Pimsner algebra, every covariant representation $(\pi , \wE_\pi)$ induces a morphism $\toe \ra L(\mM_\pi)$, with image the \sC subalgebra of $L(\mM_\pi)$ generated by $\wE_\pi$. If $\wE_\pi$ has support the identity, then we obtain a monomorphism $\coe \hra L(\mM_\pi)$ (\cite[Thm.3.12]{Pim97}).

Let us now introduce the inductive limit $\mA_\infty := \mA \stackrel{\rho}{\ra} \mA \stackrel{\rho}{\ra} \cdots$.

\begin{lem}
Let $\mE \ra X$ be a vector bundle. There exist covariant representations $(\pi , \wE_\pi)$ of $(\mA,\rho)$ with rank $\mE$ if and only if $\mA_\infty \neq \left\{ 0 \right\}$.
\end{lem}

\begin{proof}
We proceed as in \cite[Prop.2.2]{Sta93}. Let $\mA_\infty \neq \left\{ 0 \right\}$. We consider the crossed product $\cpn$ of $\mA$ by $\rho$, with the nondegenerate morphism $i_\mA : \mA \hra \cpn$. Recall that $i_\mA$ extends to the multiplier algebras; furthermore, there is a partial isometry $v \in M(\cpn)$ with support $p := i_\mA (1_{M(\mA)})$, and the relation $i_\mA \circ \rho (a) = v i_\mA (a) v^*$, $a \in \mA$, holds. Note that every $f \in C(X)$ defines a multiplier on $\cpn$, say $i_\mA (f)$. Now, since $\rho (fa) = f \rho (a)$, we find $v i_\mA (fa) = i_\mA (f) v i_\mA (a)$, so that $i_\mA (f)$ commutes with $i_\mA(a)$, $v$. Furthermore, since $C(X) \cdot \mA$ is dense in $\mA$, and $\mA_\infty \neq \left\{ 0 \right\}$, we find that $i_\mA (C(X)) \cdot \cpn$ is dense in $\cpn$. Thus, $\cpn$ is a $C(X)$-algebra. By \cite[\S 1.6]{Kas88}, there exists a faithful $C(X)$-module representation $\nu$ of $\cpn$ over a Hilbert $C(X)''$-module $\mM$. For the same reason, there is a faithful $C(X)$-module representation $\nu'$ of $\coe$ over a Hilbert $C(X)''$-module $\mM'$. We consider the bimodule tensor product $\mM \otimes_{C(X)''} \mM'$, and claim that the pair $(\pi , \wE_\pi)$, $\pi:= (\nu \circ i_\mA) \otimes 1$, $\wE_\pi := \nu(v) \otimes \nu' (\wE)$, is a covariant representation of $\mA$ over $\mM \otimes_{C(X)''} \mM'$. In fact, it is obvious that $\wE_\pi$ is isomorphic to $\wE$ as a Hilbert $C(X)$-bimodule; furthermore, if $\left\{ \psi_l \right\}$ is a finite set of generators for $\nu' (\wE)$, by definition we find that $\wE_\pi$ has support 
\[
\sum_l ( \nu(v) \otimes \psi_l) \cdot ( \nu(v)^* \otimes \psi_l^*) = (\nu (v) \cdot \nu(v^*)) \otimes 1 = \nu (p) \otimes 1 \,
\]
\noindent If $a \in \mA$, $\psi \in \nu'(\wE)$, then
\[
\begin{array}{ll}

\left( \nu \circ i_\mA (\rho(a)) \otimes 1 \right)  \cdot  \left( \nu(v) \otimes \psi \right) & = 

\nu \left(  i_\mA (\rho(a)) \cdot v   \right)  \otimes \psi  = \\ & =

\nu \left(  v \cdot i_\mA(a) \cdot v^* \cdot v   \right)  \otimes \psi   = \\ & =

\left( \nu ( v ) \otimes \psi \right)  \cdot \left( \nu \circ i_\mA(a) \otimes 1 \right) \ \ ;

\end{array}
\]
\noindent so that $\pi \circ \rho (a)  \varphi  = \varphi \pi (a)$, $\varphi \in \wE_\pi$. Since by definition $\varphi \pi (a) = \sigma_\pi \circ \pi (a) \varphi$, where $\sigma_\pi$ is the inner endomorphism induced by $\wE_\pi$, we conclude that $\pi \circ \rho = \sigma_\pi \circ \pi$, and the first implication is proven. Viceversa, let $\mA_\infty = \left\{ 0 \right\}$, so that for every $a \in \mA$ there is $k \in \bN$ such that $\rho^k(a) = 0$. Then, if $(\pi , \wE_\pi)$ is a covariant representation, we obtain $0 = \pi(\rho^k(a)) = \sum_{ \left|  \vL  \right| = k } \psi_\vL \pi (a) \psi_\vL^*$, where $\left\{ \psi_l \right\}$ is a set of generators of $\wE_\pi$ and $\psi_\vL := \psi_{l_1} \cdots \psi_{l_k}$. Note that if $d$ is the rank of $\mE$, then $\sum_{\vL = \left| k \right|} \psi_\vL^* \psi_\vL = \sum_L \left \langle \psi_\vL , \psi_\vL \right \rangle = d^k$. Now, $\pi (a) = d^{-k} \sum_L \psi_\vL^* \pi(a) \psi_\vL = 0$, thus $\pi = 0$. This is a contradiction, so that the lemma is proven.
\end{proof}

\begin{prop}
\label{prop_cp}
Let $\mA$ be a $C(X)$-algebra, $\rho$ a $C(X)$-endomorphism of $\mA$ with $\mA_\infty \neq \left\{ 0 \right\}$. Then, for every vector bundle $\mE \ra X$ there exists up to isomorphism a unique $C(X)$-algebra $\cpen$ such that

\begin{itemize}

\item there is a non-degenerate $C(X)$-morphism $i : \mA \ra \cpen$;

\item $\wE$ is contained in $\cpen$ as a finitely generated Hilbert $C(X)$-bimodule; furthermore, $i \circ \rho = \sigma_\mE \circ i$, where $\sigma_\mE$ is the inner endomorphism induced by $\wE$ on $\cpen$;

\item $\cpen$ is generated as a \sC algebra by $i(\mA)$, $\wE$;

\item for every covariant representation $(\pi , \wE_\pi)$ there is a unique non degenerate $C(X)$-representation $\Pi : \cpen \ra L(\mM_\pi)$ such that $\Pi \circ i = \pi$ and $\Pi (\wE) = \wE_\pi$.

\end{itemize}

\end{prop}

\begin{proof}
The unicity of the crossed product follows the universality property w.r.t. covariant representations. Let us now consider the *-algebra $\mB$ generated by $\mA$, $\wE$ with relations
\begin{equation}
\label{eq_uni_cp}
\left\{ 
\begin{array}{l}
\psi^* \psi' = \left\langle \psi , \psi' \right\rangle \\
f \psi = \psi f \\ 
(af) \cdot \psi = a \cdot (f \psi)  \\
\psi \cdot a = \rho (a) \cdot \psi
\end{array}
\right.
\end{equation}
\noindent where $a \in \mA$, $\psi , \psi' \in \wE$, $f \in C(X)$ and $\left\langle \cdot , \cdot \right\rangle$ is the $C(X)$-valued scalar product on $\wE$. Every covariant representation $(\pi , \wE_\pi)$ induces a representation of $\mB$ over the Hilbert $C(X)''$-module $\mM_\pi$. We denote by $\cpen$ the \sC algebra obtained by endowing $\mB$ with the maximal seminorm w.r.t. such representations. It follows from the previous lemma that the natural $C(X)$-morphism $i : \mA \ra \cpen$ is nondegenerate, $\mA_\infty$ being $\neq \left\{ 0 \right\}$. Moreover, $\wE$ is contained in $\cpen$ as a finitely generated Hilbert $C(X)$-bimodule. By the last of (\ref{eq_uni_cp}), the inner endomorphism induced by $\wE$ extends $\rho$ as desired.
\end{proof}

The next proposition states the compatibility of our crossed product with those defined by Stacey, Paschke, Cuntz, by choosing trivial vector bundles.

\begin{prop}
Let $\mA$ be a $C(X)$-algebra, $\rho$ a $C(X)$-endomorphism of $\mA$, ${\mE_d} := X \times \bC^d$ the trivial rank $d$ vector bundle over $X$. Then $\mA \rtimes_\rho^{{\mE_d}} \bN$ is isomorphic to the crossed product $\mA \rtimes_\rho^d \bN$ (in the sense of \cite{Sta93}).
\end{prop}

\begin{proof}
The module of continuous sections of ${\mE_d}$ is isomorphic to the free rank $d$ Hilbert $C(X)$-module. Thus, we can pick $d$ orthonormal generators $\psi_1 , \dots , \psi_d$, that appear in $\mA \rtimes_\rho^{{\mE_d}} \bN$ as isometries such that $\psi_i^* \psi_j = \delta_{ij}1$, and inducing the inner endomorphism extending $\rho$. Of course both $\mA \rtimes_\rho^{{\mE_d}} \bN$ and $\mA \rtimes_\rho^d \bN$ are generated by elements of the type $i(a) \psi_\vL \psi_\vM^*$; thus, we have to verify only that $\mA \rtimes_\rho^d \bN$ satisfies the universal property w.r.t. covariant representations $(\pi,\wE_{d,\pi})$ (note that $\mA \rtimes_\rho^d \bN$ is a $C(X)$-algebra, so it makes sense to consider $C(X)$-representations of $\mA \rtimes_\rho^d \bN$). Let now $\mM_\pi$ be the Hilbert $C(X)''$-module carrying the covariant representation $(\pi,\wE_{d,\pi})$; we denote by $\psi_{i , \pi}$, $i = 1 , \ldots , d$, the orthogonal isometries in $L(\mM_\pi)$ generating $\wE_{d,\pi}$. We have to prove that there exists a $C(X)$-module representation $\Pi$ of $\mA \rtimes_\rho^d \bN$ as in Prop.\ref{prop_cp}. We consider a faithful state $\lambda$ on $C(X)''$, and the Hilbert space $\mM_{\pi,\lambda}$ obtained by introducing on $\mM_\pi$ the $\bC$-valued scalar product $(v,w)_\lambda := \lambda( v,w )$, $v,w \in \mM_\pi$. Thus, there is a unital monomorphism $i_\lambda : L(\mM_\pi) \hra L(\mM_{\pi,\lambda})$, and $(i_\lambda \circ \pi , \left\{ i_\lambda(\psi_{i , \pi}) \right\}_{i=1}^d )$ is a covariant representation in the sense of \cite[\S 2]{Sta93}. By construction of $\mA \rtimes_\rho^d \bN$, there exists a representation $\Pi_\lambda$ of $\mA \rtimes_\rho^d \bN$ such that $i_\lambda \circ \pi  = \Pi_\lambda \circ i$ and $\Pi_\lambda (\psi_i) = i_\lambda (\psi_{i , \pi})$. Now, since $i_\lambda$ is injective and $\Pi_\lambda (\mA \rtimes_\rho^d \bN) \subseteq i_\lambda(L(\mM_\pi))$, we can define the \sC algebra morphism $\Pi := \left. i_\lambda^{-1} \right|_{i_\lambda (L(\mM_\pi))} \circ \Pi_\lambda$. If $f \in C(X)$, $b \in \cpen$, then $\Pi (i(f)b) = i_\lambda^{-1} \circ \Pi_\lambda (i(f)) \ \Pi (b) = \pi (f) \Pi (b)$; since $\pi$ is a $C(X)$-representation, we conclude that $\Pi$ is a $C(X)$-representation. It is now clear that $\Pi$ satisfies by construction the required properties, thus the proposition is proven.
\end{proof}

\begin{rem} {\em
Let $\rho$ be an automorphism of a \sC algebra $\mA$, $\mL \ra \spzro$ a line bundle. Then, the inner endomorphism $\sigma_\mL$ is a 'locally unitary' automorphism of $\mA \rtimes_\rho^\mL \bN$, in analogy with the notion introduced in \cite{PR84}.
} \end{rem}

\begin{ex} {\em
\label{ex_ct}
Let $\mA$ be a stable continuous trace \sC algebra with compact spectrum $X$. Then, for every $C(X)$-endomorphism $\rho$ of $\mA$ there exists a unique (up to isomorphism) finitely generated Hilbert $C(X)$-bimodule $\wE (\rho)$ contained in the multiplier algebra $ M (\mA)$, inducing $\rho$ on $\mA$ in the usual way (\ref{12}) (see \cite[\S 4]{Hir01}). $\wE (\rho)$ is isomorphic to the module of continuous sections of a vector bundle $\mE (\rho) \ra X$. The \sC subalgebra of $M (\mA)$ generated by $\mA$ and $\wE (\rho)$ is isomorphic to the crossed product $\mA \rtimes_\rho^{\mE (\rho)} \bN$, since it obviously satisfies the universal property w.r.t. covariant representations.
} \end{ex}

\begin{ex} {\em
\label{ex_cp}
Let $\mE \ra X$ be a vector bundle, and $\coe^0$ the fixed point \sC algebra of $\coe$ w.r.t. the circle action $\wa z (\psi):= z \psi$, $\psi \in \wE$, $z \in \bT$. By general facts (see for example \cite[\S 3]{DPZ97}), $\coe^0$ can be regarded as the inductive limit $\lim \limits_{\ra r} (\mE^r , \mE^r)$, where $(\mE^r , \mE^r)$ is embedded in $(\mE^{r+1} , \mE^{r+1})$ by tensoring on the right by the identity of $(\mE , \mE)$. If $p \in (\mE , \mE)$ is a projection, we define the shift endomorphism $\wa p := t \mapsto p \otimes t$, with $t \in (\mE^r,\mE^r) \subset \coe^0$. Let in particular $p$ be a rank $1$ projection, and $\rho := \wa p$. We prove that 
\[ \coe \simeq \coe^0 \rtimes_\rho^\mL \bN \ , \]
\noindent where $\mL := p \mE \subset \mE$ is the line bundle projected by $p$. Since 
\begin{equation}
\label{eq_we}
\wE = (\mE , \mE) \cdot \wa \mL 
:= {\mathrm {span}} \left\{ \right.
t \varphi \ , \ t \in (\mE , \mE) , \varphi \in \wa \mL
\left. \right\}
\ ,
\end{equation}
\noindent we find that $\coe$ is generated by $\coe^0$, $\wa \mL$. Furthermore, $\rho$ becomes inner in $\coe$: in fact, the identity $(p \otimes t) \cdot (\varphi \otimes^r 1) = \varphi \otimes t$, $t \in (\mE^r,\mE^r)$, $\varphi \in p \wE \simeq \wL$, is equivalent to $\rho (t) \varphi = \varphi t$ in the CP-algebra $\coe$. Finally, we have to prove that $\coe$ satisfies the universal property for covariant representations. Let $(\pi : \coe^0 \ra L(\mM_\pi) , \wL_\pi)$ be a covariant representation. Since $\rho$ is induced by $\wa \mL$ we find $(\iota , \rho) = (\iota , \mL)$, so that $(\rr) = (\mL , \mL)$; since $p$ has rank $1$, we find $(\mL , \mL) = p \cdot (\mE ,\mE) \cdot p = C(X) \cdot p$. If $t,t' \in (\mE , \mE)$, $\varphi , \varphi' \in \wL$, we find, with $f p := pt^* t'p \in C(X) \cdot p$,
\[
\varphi^* t^* t' \varphi' = 
\varphi^* p t^* t' p \varphi' = 
f \varphi^* \varphi' \in C(X) \ ;
\] 
\noindent moreover, if $\varphi_\pi , \varphi'_\pi \in \wL_\pi$ are the images of $\varphi , \varphi' \in \wL$ w.r.t. the isomorphism $\wL \simeq \wL_\pi$,
\begin{equation}
\label{eq_ep}
\varphi^*_\pi \pi(t^* t') \varphi'_\pi = 
\varphi^*_\pi \pi(p) \pi(t^* t') \pi(p) \varphi'_\pi = 
f \varphi^*_\pi \varphi'_\pi = f \varphi^* \varphi'  \in C(X) \ .
\end{equation}
\noindent The previous equality implies that $\wE_\pi := \pi (\mE , \mE) \cdot \wL_\pi$ is a Hilbert $C(X)$-bimodule in $L(\mM_\pi)$. Recall from (\ref{eq_we}) that if $\psi \in \wE$, then $\psi = t \varphi$, $t \in (\mE , \mE)$, $\varphi \in \wL$ (up to linear combinations with coefficients in $C(X)$). We define the following morphism of Hilbert $C(X)$-bimodules:
\[
\phi : \wE \ra \wE_\pi \ \ : \ \ t \varphi \mapsto \pi (t) \varphi_\pi \ \ .
\]
\noindent (\ref{eq_ep}) implies that $\phi$ preserves the $C(X)$-valued scalar product, thus $\phi$ is an isomorphism. By universality of the CP-algebra, there exists a unique $C(X)$-monomorphism $\Pi : \coe \hra L(\mM_\pi)$ extending $\phi$: by construction, $\Pi$ is the desired representation extending $\pi$.
} \end{ex}

\section{The \sC algebra of a $G$-vector bundle.}
\label{cavb_gb}

We start the present section by fixing some notations: if $d \in \bN$, then $H \simeq \bC^d$ is the standard Hilbert space with rank $d$; $\bM_d$ is the \sC algebra of $d \times d$ matrices with coefficients in $\bC$; $\bM_{d,d'}$, $d' \in \bN$, is the Banach space of $d \times d'$ matrices with coefficients in $\bC$. If $r \in \bN$, we denote by $H^r$ the $r$-fold tensor power of $H$; if $r = 0$ we define $H^0 := \bC$. We also make use of the notation $(\hrs)$, $r,s \in \bN$, to denote the Banach space of linear operators from $H^r$ into $H^s$, so that $(\hrs) \simeq \bM_{d^r , d^s}$. Note that the Cuntz algebra $\mO_d$ is generated by the spaces $(\hrs)$, $r,s \in \bN$, according to the procedure described in \S \ref{preli}.

For basic notions and terminology about vector bundles, we refere to \cite{Ati,Kar}, while for generic fibre bundles we refere to \cite{Hus}.

Let $X$ be a compact Hausdorff space, $\mE \ra X$ a rank $d$ vector bundle; we denote by $\mE_x \simeq H$ the fibre of $\mE$ over $x \in X$, and by $\wE$ the Hilbert $C(X)$-bimodule of continuous sections of $\mE$. We denote the elements of $\wE$ by
\begin{equation}
\label{def_px}
\psi : X \ra \mE \ \ : \ \ 
x \mapsto \psi_x \in \mE_x \subseteq \mE  \ \ .
\end{equation}
\noindent In the sequel, we will often make use of the following property: let $U \subseteq X$ be a closed set with nonempty interior, $\pi_U : \mE |_U \ra U \times H$ a local chart. Then, by functoriality a local chart $\pi^r_i : \mE^r |_{U} \ra U \times H^r$ is induced for every $r \in \bN$. Now, for every $r \in \bN$ the bimodule of continuous sections of $U \times H^r$ is given by the free bimodule $C(U) \otimes H^r$, that we identify with the space of continuous maps from $U$ into $H^r$. The local charts $\pi_U^r$ define restriction morphisms
\[
\pi_U^{0,r} : (\iota , \mE^r) \ra C(U) \otimes H^r    \ \ , \ \ 
\pi_U^{0,r} (\psi) := 
\left\{ x \mapsto  \pi_U^r \circ \psi_x \right\}
\ ,
\]
\noindent $\psi \in \wa {\mE^r} \simeq (\iota , \mE^r)$, $x \in U$. More in general, by extending over elements of the type $\psi_\vL \psi_\vM^*$, $\left| L \right| = s$, $\left| M \right| = r$ (see \S \ref{preli}), we obtain restriction morphisms
\begin{equation}
\label{eq_rm}
\pi_U^{r,s} : (\ers) \ra C(U) \otimes (\hrs) \ \ ,
\end{equation}
\noindent Note that $\pi_U^{r,q} (t't) =  \pi_U^{s,q} (t') \ \pi_U^{r,s} (t)$, $t \in (\ers)$, $t' \in (\mE^s , \mE^q)$. Thus, the following local chart is induced on $\coe$:
\begin{equation}
\label{def_flc}
\alpha_U : \coe |_U \ra C(U) \otimes \mO_d \ \ , \ \
\alpha_U (t) := \pi_U^{r,s} (t) \ \ ,
\end{equation}
\noindent $t \in (\ers) \subset \coe$.

\subsection{Unitary operators on vector bundles.}
Let $\mUE$ be the group of unitary endomorphisms of $\mE$. It is well-known that $(\mE,\mE)$ is a locally trivial continuous bundle with fibre $\bM_d$, thus $\mUE$ can be naturally regarded as the group of continuous sections of a fibre bundle $\mcUE \ra X$, that we call {\em the unitary bundle of} $\mE$. $\mcUE$ can be constructed in the following way: let $\left\{ u_{ij} : X_i \cap X_j \ra \ud \right\} \in H^1 ( X , \ud)$ be a set of {\em transition maps} associated with $\mE$ for an open trivializing cover $\left\{ X_i \right\}$ (see \cite[I.3.5]{Kar} about such a terminology); then, $\mcUE$ is constructed by clutching the bundles $X_i \times \ud$ {\em via} the maps $f_{ij} (x , u) := (x , u_{ij}(x) \cdot u \cdot u_{ij}^* (x) )$, $(x , u) \in (X_i \cap X_j) \times \ud$. We denote by $\eta : \mcUE \ra X$ the surjective map associated with the bundle structure of $\mcUE$. By construction, $\mcUE$ has fibre $\mcUE_x := \eta^{-1}(x) \simeq \ud$ and structure group $\left. \ud \right/ \bT$, acting by adjoint action. By functoriality, if $\pi_U : \mE |_U \ra U \times H$ is a local chart for $\mE$, then a local chart 
\begin{equation}
\label{def_lcud}
\eta_U : \mcUE |_U \ra U \times \ud  \ \
\end{equation}
\noindent is induced. The elements of $\mUE$ can be described by families of continuous maps $g_i : X_i \ra \ud$, satisfying the cocycle relations $u_{ij} g_j = g_i u_{ij}$. If $\mE$ is trivial, then $\mcUE \simeq X \times \ud$ and $\mUE$ is isomorphic to $C(X , \ud)$, i.e. the group of continuous maps from $X$ into $\ud$.

\begin{rem} {\em
Let $U \subseteq X$ be a closed set with non empty interior trivializing $\mE$, $\pi_U^{1,1} : (\mE , \mE ) \ra C(U) \otimes (H,H)$ be the restriction morphism (\ref{eq_rm}). In particular, we obtain the restriction morphism 
\[  
\pi_U^{1,1} : \mUE \ra C(U , \ud ) \ \ .
\]
} \end{rem}

\begin{rem} {\em
\label{rem_full}
We denote by 
\[ g_x \in \mcUE_x \simeq \ud \]
\noindent the evaluation of $g \in \mUE$ over $x \in X$. It is well-known that $\mcUE$ is {\bf full}, i.e. for every $x \in X$, $u \in \mcUE_x$ there is $g \in \mUE$ such that $u = g_x$. The argument for the proof is the following. Pick an open set $U \subseteq X$ with a local chart $\eta_U : \mcUE |_U \ra U \times \ud$; if $(\omega , g_0) \in U \times \ud$, then there is a continuous map $u : [0,1] \ra \ud$ such that $u(1) = g_0$, $u(0) = 1$. We consider a cutoff $\lambda \in C(X)$, $0 \leq \lambda \leq 1$, with support contained in $U$, $\lambda (\omega) = 1$, and construct a continuous section $g \in \mUE$, 
\[
g_x :=
\left\{
\begin{array}{ll}
\eta_U^{-1} ( x \ ,\ u \circ \lambda (x) )  \ \ , \ \   x \in U   \\
1 \ \ , \ \ x \in X - U     \\
\end{array}
\right.
\]
\noindent where $1 \in \mUE$ is the identity. Note that $g_{\omega} = \eta_U^{-1} (\omega , g_0 )$. Since $g_0$ is an arbitrary element of $\ud \simeq \mcUE_{\omega}$, the assertion is proved.
} \end{rem}

\bigskip

We consider the {\it determinant map} $\det : \mUE \ra C(X,\bT)$, $(\det(g)) (x) := \det(g_x)$, and introduce the group of special unitaries $\mSUE := \det^{-1} \left\{ 1 \right\}$. It is well-known that $\mSUE$ is the group of continuous sections of a bundle $\mcSUE \ra X$ with fibre $\sud$.

Let $\mG$ be a closed subgroup of $\mUE$; then $\mG$ acts in a natural way over $\mE$, which becomes a $\mG$-vector bundle (with trivial $\mG$-action over $X$) in the sense of \cite[\S 1.6]{Ati}. Note that in general $\mG$ is not (locally) compact (for example, consider $\mE := X \times H$, so that $\mUE = C(X,\ud)$). A technical consequence of this fact is that there is no immediate analogue of the Haar measure, for such a kind of groups.

\begin{rem} {\em
\label{rem_glc}
Let $G_0$ be a locally compact group, $\mE \ra X$ a $G_0$-vector bundle such that the $G_0$-action is trivial on $X$. Then, every $g \in G_0$ defines a unitary map on $\mE$, so that there is a morphism $U : G_0 \ra \mUE$ of topological groups. If $U$ is injective, then $UG_0 \subset \mUE$ is locally compact.
} \end{rem}

Let $\mG \subseteq \mUE$ be a closed group; by the identification $(\mE , \mE) \simeq (\wE , \wE)$, we obtain that $\mG$ acts on $\wE$ by unitary $C(X)$-bimodule operators.

Let $\coe$ be the CP-algebra associated with $\wE$. We introduce a canonical $\mUE$-action on $\coe$: if $g \in \mUE$, then $g^{{\otimes}^r} \in (\mE^r , \mE^r)$ for every $r \in \bN$, and we define
\begin{equation}
\label{24} 
\wa{g} (t) := g^{{\otimes}^s} \cdot t \cdot {g^*}^{{\otimes}^r}  \in (\ers) \ ,
\end{equation}
\noindent where $g \in \mUE$, $t \in ({\mE}^r ,{\mE}^s)$. The map $\left\{ g \mapsto \wa g \right\}$ defines an $\mUE$-action by $C(X)$-automorphisms on $\coe$ (this fact is a consequence of the universality of the CP-algebra, see \cite[\S 3]{DPZ97}). We also consider the fixed-point Banach $C(X)$-bimodules
\begin{equation}
\label{def_ersg}
(\ers)_\mG := \left\{y \in (\ers) : \wa g (y) = y , g \in \mG \right\} \  , \ r,s \in \bN \ .
\end{equation}
\noindent If $X$ reduces to a single point, then we recover the canonical $\ud$-action over the Cuntz algebra $\mO_d$ studied in \cite[\S 1]{DR87}; if $G_0 \subseteq \ud$ is a closed group we denote by $(\hrs)_{G_0}$, $r,s \in \bN$, the corresponding fixed-point Banach spaces, and by $\mO_{G_0} \subseteq \mO_d$ the associated \sC algebra.

Let $\left\{ t_x \in \mO_d \right\}_{x \in X}$ be the vector field associated with $t \in \coe$. Then, it is clear that the vector field associated with $\wa g (t)$ is given by $\left\{ \wa g_x (t_x) \in \mO_d \right\}_{x \in X}$.

As we will see in the sequel, the Banach $C(X)$-bimodules $(\ers)_\mG$ are not necessarily finitely generated, thus they do not correspond to modules of continuous sections of vector bundles in the usual sense. They define more in general (non-locally trivial) {\em Banach bundles} in the sense of \cite{Dup74}, also called {\em quasi vector bundles} in \cite[I.1]{Kar}. In the special case in which $G_0$ is a locally compact group and $\mE$ is a $G_0$-vector bundle (with trivial $G_0$-action over $X$, see Rem.\ref{rem_glc}), then the bimodules $(\ers)_{G_0}$ correspond to vector bundles (as can be proven with the argument used in \cite[Prop.1.6.2]{Ati}).

\subsection{Algebraic properties of $\cog$.}
We consider the *-subalgebra $\aog$ of $\coe$ generated by the Banach $C(X)$-bimodules $(\ers)_\mG$, $r,s \in \bN$, and denote by $\cog$ the closure of $\aog$ in $\coe$. $\cog$ naturally inherits from $\coe$ the structure of continuous bundle of \sC algebras, with fibre 
\[
(\cog)_x \subseteq \mO_d \ \ .
\]
\noindent If $y \in \cog$, we denote by $y_x \in \mO_d$ the evaluation of $y$ over $x \in X$ as a vector field. As by (\ref{def_sigma}), we denote by $\sigma$ the canonical endomorphism of $\coe$. If $g$ is any element of $\mUE$, then for $t \in (\ers)$ we find $\wa g \circ \sigma (t) = \wa g (1 \otimes t) = 1 \otimes \wa g (t)$, so that $\wa g \circ \sigma = \sigma \circ \wa g$. Thus, for every $\mG \subseteq \mUE$, the restriction $\sigma_\mG \in {\bf end} \cog$ of the canonical endomorphism is well defined. For intertwiners, as usual we use the notation
\[
(\sigma_\mG^r,\sigma_\mG^s) := \left\{ y \in \cog : y \sigma_\mG^r(y') = \sigma_\mG^s(y') y , y' \in \cog  \right\} .
\]

\bigskip

We now discuss structural properties of $\cog$ and the canonical endomorphism $\sigma_\mG$. Let us denote by $\lambda \mE \subseteq \mE^d$ the exterior tensor product with order the rank of $\mE$. $\lambda \mE \ra X$ is a line bundle, and $\lambda (\mE \oplus \mE') = \lambda \mE \otimes \lambda \mE'$ (see for example \cite[\S 9.12 (b)]{Kar} about the above identity); by identifying $H^2(X, \bZ)$ with the group of isomorphism classes of line bundles over $X$ (endowed with the tensor product as group operation), we obtain an epimorphism $\lambda : K(X) \ra H^2(X,\bZ)$, i.e. the first Chern class. We now give a local description of $\lambda \mE$. We consider a trivialization $\pi_i : \left. {\mE} \right|_{X_i} \mapsto X_i \times H$, such that $\mE$ admits a family of transition maps $\left\{ u_{ij} := \pi_i \circ \pi_j^{-1} : X_i \cap X_j \ra \ud \right\}$; note that for every $r \in \bN$, a local chart $\pi^r_i : \mE^r |_{X_i} \ra X_i \times H^r$ is induced. We pick a partition of unity $\left\{ \lambda_i \right\}$ subordinate to the open cover $\left\{ X_i \right\}$ and consider, for each index $i$, the continuous sections 
\begin{equation}
\label{eq_dr}
R_i := \lambda_i \cdot (\pi^d_i)^{-1} \circ R \in (\iota , \mE^d) \ ,
\end{equation}
\noindent where
\begin{equation}
\label {21} 
R := \sum\limits_{p \in \bP_d} \mathrm{sign} \left( p \right) e_{p\left( 1\right) }\otimes \cdots \otimes e_{p\left( d\right) } \ \in (\iota , H^d)
\end{equation}
\noindent is regarded as a constant section of $X_i \times H^d$ ($\bP_d$ denotes the permutation group, and $\left\{ e_h \right\}$ is the standard orthonormal basis of $H$). The $R_i$'s satisfy the relations
\begin{equation}
\label {221} 
\begin{array}{ll}
R_i^* R_j =
\left\langle R_i , R_j \right\rangle = 
      \lambda_i \lambda_j \left\langle R , u_{ij}^{\otimes^d} R \right\rangle =
      \lambda_i \lambda_j \det \left( u_{ij} \right) \in C(X) \ , \\
\end{array}
\end{equation}
\begin{equation}
\label {222} \lambda_i R_j = \lambda_j \det \left( u_{ij} \right) R_i \ ,
\end{equation}
\begin{equation}
\label{223}
\sum_i R_i^* R_i = \sum_i \left \langle R_i , R_i \right \rangle = 1 \ ,
\end{equation}

\noindent and generate $(\iota , \lambda \mE) \subset (\iota , \mE^d)$ as a right Hilbert $C(X)$-module. Since $\wa g (R_i) = \det (g) R_i = R_i$, $g \in \mSUE$, the $\mSUE$-action on $(\iota , \lambda \mE)$ is trivial. Now, the support of $(\iota , \lambda \mE)$ in $\coe$ is the projection
\begin{equation}
\label {23} P := \sum \limits_i R_i R_i^{*} \ ;
\end{equation}
\noindent so that, by using (\ref {23}), we find $\wa g (P) = \det (g) P \det(g^*) = P$ for $g \in \mUE$. The previous considerations imply that $ (\iota , \lambda \mE) \subset \mO_{\mSUE} $ and $P \in \mO_{\mUE} $.

An element of $({\mE}^2 , {\mE}^2)$ that will play a special role in the sequel is the symmetry $\theta \in (\mE^2,\mE^2)$,
\begin{equation}
\label{defexc}
\theta (\psi \otimes \psi') := \psi' \otimes \psi \ ,
\end{equation}
\noindent where $\psi , \psi' \in (\iota,\mE)$. It is clear that $\theta = \theta^* = \theta^{-1}$. As an element of $\coe$, $\theta$ can be expressed as
\begin{equation}
\label{defsim}
\theta = \sum \limits_{l,m} {\psi}_m {\psi}_l {\psi}^*_m {\psi}^*_l \ ,
\end{equation}
\noindent where $\left\{ \psi_l \right\}$ is a set of generators of $(\iota,\mE)$. Since the right $C(X)$-action coincides with the left one on $(\iota,\mE)$, we obtain that $\theta$ is actually well defined (i.e., it does not depend on the choice of generators), bypassing the problems mentioned in \cite[Rem.6.4]{BL97}. Thus, for $g \in \mUE$ we find
\[
\wa g (\theta) = 
\sum \limits_{l,m} (g{\psi}_m) (g{\psi}_l) (g{\psi}_m)^* (g{\psi}_l)^* =
\theta \ ,
\]
\noindent so that $\theta \in \mO_{\mUE}$.

Let $\bP_\infty$ denote the group of finite permutations of $\bN$, $\bS \in {\bf end} \bP_\infty$ the shift endomorphism $(\bS p)(1) := 1$, $(\bS p)(n) := 1 + p(n-1)$, $p \in \bP_\infty$. By using the canonical endomorphism and the symmetry we get, with the same methods used in \cite[\S 2]{DR87}, the following representation of $\bP_\infty$ in $\coe$:
\begin{equation}
\label{def_thetap}
\bP_r \ni p 
\mapsto 
\theta (p) \in (\mE^r , \mE^r) : \; \; 
\theta (p) \cdot \otimes_i^r \psi_i := \otimes_i^r \psi_{p(i)} \ , \; \; \; 
\psi_i \in (\iota , \mE) \ .
\end{equation}
\noindent Note that $\theta \circ \bS(p)= \sigma \circ \theta (p)$, and that each $\theta (p)$ is a word in $\theta$, $\sigma (\theta)$, $\cdots$, $\sigma^s (\theta)$ for some $s \in \bN$. Since the canonical endomorphism commutes with the $\mUE$-action, $\theta (p)$ belongs to $\mO_{\mUE}$ for each $p$. In particular, we will make use of the unitary operators $\theta (r,s) \in ( \mE^{r+s} , \mE^{r+s} )$, permuting the first $r$ factors of the tensor product $\mE^{r+s}$ with the remaining $s$. Elementary properties of these operators, which are the same for vector bundles as well as for Hilbert spaces, are given in \cite[\S 2]{DR87}. We will make use of the obvious identity
\[
\theta (s,1) \cdot (t \otimes 1) = (1 \otimes t) \cdot \theta (r,1) \quad ,
\]
\noindent $t \in (\ers)$, which regarded in $\coe$ becomes
\begin{equation}
\label{ps_theta}
\theta (s,1) \cdot t = \sigma(t) \cdot \theta (r,1) \ .
\end{equation}
\noindent  Note that in particular $\theta \psi = \sigma (\psi)$, $\psi \in (\iota,\mE)$.

\begin{rem} {\em
Let $R \in (\iota , H^d)$ be the isometry defined in (\ref{21}). Then, the identity $RR^* = \frac 1{d!} \sum \limits_{p \in \bP_d} \mathrm{sign} (p) \theta (p)$ holds (see remarks after \cite[Lemma 3.7]{DR87}); by applying (\ref{23}), we find
\begin{equation}
\label{27}
P = \sum \limits_i R_i R_i^*
  = \frac 1{d!} \sum \limits_{p \in \bP_d} \mathrm{sign} (p) \theta (p) \ ,
\end{equation}
\noindent so that $P$ belongs to the $\sigma$-stable algebra generated by $\theta$.
} \end{rem}

\begin{prop}
\label{rel_com_og}
\
\begin{itemize}

\item $\sigma (t) = \lim_{r \ra \infty} \theta (r+k,1) \cdot t \cdot \theta (r,1)^*$, for every $t \in \mO_\mE^k$;

\item if $t \in \coe$ and $\left[ t , \theta (p) \right] = 0$ for $p \in \bP_\infty$, then $\sigma (t) = t$ and $t \in C(X)$;

\item if $\mG$ is any subgroup of $\mUE$, then $\cog' \cap \coe = C(X)$;

\item the amenability property holds, i.e. $(\ers)_\mG = (\sigma_\mG^r,\sigma_\mG^s)$ for $r,s \in \bN$.

\end{itemize}
\end{prop}

\begin{proof}
\
\begin{itemize}

\item It follows immediately by (\ref{ps_theta}), and by considering the inductive limit.

\item If $t$ commutes with $\theta (p)$ for every $p \in \bP_\infty$, then $t_x \theta_x(p) = \theta_x(p) t_x$ for each fibre $t_x \in \mO_d$ and the corresponding permutation operator $\theta_x(p) \in \mO_d$, as $x$ varies in $X$. Thus, by \cite[Lemma 3.2]{DR87} we find $t_x \in \bC$, and $t \in C(X)$.

\item Since $\theta \in \mO_{\mUE} \subset \cog$ and $\cog$ is $\sigma_\mG$-stable, we find that $t \in \cog' \cap \coe$ commutes with $\theta (p)$ for every $p \in \bP_\infty$, so that $t \in C(X)$.

\item If $t \in (\sigma_\mG^r,\sigma_\mG^s)$, then $t_{\vL \vM} := \psi_\vL^* t \psi_\vM \in \cog' \cap \coe = C(X)$ for $\psi_\vM \in (\iota,\mE^r)$, $\psi_\vM \in (\iota,\mE^s)$. Thus, $t = \sum_{\vL \vM} t_{\vL \vM} \psi_\vL \psi_\vM^* \in (\ers)$, and it is by definition $\mG$-invariant. Viceversa, if $t \in (\ers)_\mG$, then $t \in (\ers) = (\sigma^r,\sigma^s)$ (see Ex.\ref{ex_ivb}); since by hypothesis $t \in \cog$, in particular $t$ belongs to $(\sigma_\mG^r,\sigma_\mG^s)$.

\end{itemize}
\end{proof}

\bigskip

\begin{rem} {\em
\label{rem_ps}
By (\ref{ps_theta}) and the amenability property, we conclude that 
\[
\theta (s,1) \cdot t = \sigma(t) \cdot \theta (r,1) \ \ , \ \ 
t \in (\sigma^r , \sigma^s) \ ;
\]
\noindent thus $\sigma$ has permutation symmetry in the sense of \cite[\S 4]{DR89A}, as for the case in which $X$ reduces to a point studied in \cite{DR87}. The same is true for every restriction $\sigma_\mG$, $\mG \subseteq \mUE$. 
} \end{rem}

\begin{lem}
\label{lem_scp}
Let $\mE \ra X$ be a rank $d$ vector bundle. Then, for every $R, R' \in (\iota , \lambda \mE)$ the following equality holds:
\begin{equation}
\label{eq_scp1}
R^* \sigma (R') = (-1)^{d-1} d^{-1} R^* R' \ .
\end{equation}
\noindent If $\mG \subseteq \mSUE$, then $\sigma_\mG \in {\bf end} \cog$ satisfies the special conjugate property in the sense of \cite[\S 4]{DR89A} if and only if the first Chern class of $\mE$ vanishes, i.e. $\lambda \mE \simeq X \times \bC$.
\end{lem}

\begin{proof}
By \cite[Lemma 2.2]{DR87}, we have 
\begin{equation}
\label{eq_scp}
R^* \sigma (R) = (-1)^{d-1} d^{-1} 1 \ ,
\end{equation}
\noindent where $R$ is defined by (\ref{21}). Let $\left\{ R_i \right\}$ be the set of generators of $(\iota , \lambda \mE)$ introduced in (\ref{eq_dr}); then $R_i^* \sigma (R_i) = (-1)^{d-1} d^{-1} \lambda_i^2$, so that, by using (\ref{221}, \ref{222}) we find $\lambda_i \lambda_j R_i^* \sigma (R_j) = (-1)^{d-1} d^{-1} \lambda_i \lambda_j R_i^* R_j$. Since $R_i$ (resp. $R_j$) has the same support of $\lambda_i$ (resp. $\lambda_j$), we can divide the previous equality by $\lambda_i \lambda_j$, and obtain
\[
R_i^* \sigma (R_j) = (-1)^{d-1} d^{-1} R_i^* R_j \ ;
\]
\noindent since $\left\{ R_i \right\}$ generates $(\iota , \lambda \mE)$, we obtain (\ref{eq_scp1}).

About the second assertion, note that $\mG \subseteq \mSUE$ implies $(\iota , \lambda \mE) \subseteq (\iota , \mE^d)_\mG$; thus, by amenability $(\iota , \lambda \mE) \subseteq (\iota , \sigma_\mG^d)$. Recall that by definition $\sigma_\mG$ satisfies the special conjugate property if and only if there is $S \in (\iota , \sigma_\mG^d)$ generating $(\iota , \lambda \mE)$ as a Hilbert $C(X)$-module, i.e., $S^* S = 1$, $SS^* = P$ (by (\ref{eq_scp1}), this implies that $S$ satisfies also (\ref{eq_scp})). This happens if and only if $\lambda \mE$ is trivial.
\end{proof}

\bigskip

The previously exposed properties (in particular the amenability) allow to give the following 

\begin{defn}
\label{def_g_dual}
Let $\mE \ra X$ be a vector bundle, $\mG \subseteq \mUE$ a closed group. We denote by $\wa \mG$ the symmetric tensor \sC category with objects the tensor powers $\mE^r$, $r \in \bN$, and arrows the invariant $C(X)$-bimodules $(\ers)_\mG$. $\wa \mG$ is called the {\bf dual of} $\mG$.
\end{defn}

Note that when $X$ reduces to a single point, then $\wa \mG$ is the category of tensor powers of the defining representation of $\mG$.

\begin{cor}
\label{cor_itc}
Let $\mG \subseteq \mUE$ be a closed group. Then, there is an isomorphism of tensor \sC categories $\wa \sigma_\mG \simeq \wa \mG$.
\end{cor}

\begin{proof}
The amenability property proved in Prop.\ref{rel_com_og} implies that $\wa \sigma_\mG$ and $\wa \mG$ are isomorphic as \sC categories. Let now $1_r \in (\mE^r , \mE^r)$ be the identity map (note that $1_r$ coincides with the identity in the Pimsner algebra $\coe$). Then, the identity
\[
t \otimes t' = (t \otimes 1_{s'}) \cdot (1_r \otimes t) = t \sigma^r (t') \ \ , \ \ 
t \in (\ers) \ , \ t' \in (\mE^{r'} , \mE^{s'} )
\]
\noindent implies that $\wa \sigma_\mG$ and $\wa \mG$ are isomorphic as tensor \sC categories.
\end{proof}

\subsection{$\cog$ as a continuous bundle.}
It is expectable that if $\mG \subseteq \mUE$ is a closed group, then a natural notion of 'fibre' of $\mG$ over $x \in X$ may be given. We now define two natural candidates for this role. As first, for every $x \in X$ we introduce the closed group
\begin{equation}
\label{def_evg}
{\mG_x} := \left\{ g_x : g \in \mG \right\} \subseteq \ud \ .
\end{equation}
\noindent Now, for every $x \in X$ there is an inclusion $(\cog)_x \subseteq \mO_d$. We consider the group
\[
{\bf aut}_{(\cog)_x} \mO_d := 
\left\{ \alpha \in {\bf aut} \mO_d  : \alpha (y) = y , y \in (\cog)_x  \right\} \ .
\]
\noindent It is clear from the definition of the canonical action (\ref{24}) that there is an immersion ${\mG_x} \hra {\bf aut}_{(\cog)_x}$: $g_x \mapsto \wa g_x$, $x \in X$. In general such an immersion is not surjective, as evident from the following example.

\begin{ex} {\em
\label{ex_point}
Let $X := [0,1]$, $\mE := X \times H$, $0 < \omega < 1$, and
\[
\mG := \left\{ g \in C(X, \ud) : g_\omega = {1} \right\} \ .
\]
\noindent Then ${\mG_x} = \ud$ for every $x \neq \omega$, and $\mG_{\omega} = \left\{ {1} \right\}$. It is clear that $(\mE^r,\mE^s) = C(X) \otimes (\hrs)$. Let now $t \in (\mE^r,\mE^s)$ be $\mG$-invariant. Then, for every neighbourhood $U$ of $\omega$, we find that $\left. t \right|_{X - U}$ is a norm-continuous map taking values into $(\hrs)_{\ud}$. Let $x \in X - U$; by \cite[Lemma 3.6]{DR87}, for $r \neq s$ we find $t_x = 0$, and for $r = s$ we find $t_x \in (H^r , H^r)_{\ud}$. By continuity, $t_\omega = 0$ if $r \neq s$, and $t_\omega \in (H^r , H^r)_{\ud}$ for $r = s$. Thus, $(\mE^r,\mE^s)_\mG = C(X) \otimes (\hrs)_{\ud}$, and $(\cog)_x = \mO_{\ud}$ for every $x \in X$. By applying \cite[Cor.3.3]{DR87}, we conclude that ${\bf aut}_{\mO_{\ud}} \mO_d \simeq \ud$ for every $x \in X$.
} \end{ex}

\begin{rem} {\em
\label{rem_dualpoint}
We consider the Cuntz algebra $\mO_d$, with the fixed-point \sC algebra $\mO_{\sud} \subset \mO_d$ w.r.t. the $\sud$-action (\ref{24}). Let $\mA$ be a $\sigma$-stable, unital \sC algebra with inclusions $\mO_{\sud} \subseteq \mA \subseteq \mO_d$ (this implies $\mA' \cap \mA \subseteq \mA' \cap \mO_d = \bC 1$, see \cite[Lemma 3.2]{DR87}). We define $\rho := \sigma |_{\mA}$, and note that $\rho$ satisfies permutation symmetry and special conjugate property in the sense of \cite[\S 4]{DR89A}. Suppose that $\mA$ is generated by the intertwiners spaces $(\rhors)$, $r,s \in \bN$. Let $\mG := {\bf aut}_\mA \mO_d$ denote the stabilizer of $\mA$ in $\mO_d$. Then, by \cite[Thm.4.1,Lemma 4.6]{DR89A}, we obtain that $\mG$ can be regarded as a closed subgroup of $\sud$ acting on $\mO_d$ by the action (\ref{24}), with $\mA = \cog$, $\rho = \sigma_\mG$, $(\rhors) = (\sigma_\mG^r , \sigma_\mG^s)$, $r,s \in \bN$.
} \end{rem}

\bigskip

Let $\mE \ra X$ be a rank $d$ vector bundle, $\mF \subseteq \coe$ a $\sigma$-stable continuous bundle of \sC algebras. We denote by $\mF_x \subseteq \mO_d$, $x \in X$, the fibre of $\mF$ over $x$, and $\rho := \sigma |_\mF \in {\bf end}\mF$.

\begin{lem}
\label{lem_ff}
With the above notation, suppose $\mO_\mUE \subseteq \mF$. Then, for every $x \in X$ the action (\ref{24}) defines a group isomorphism
\[
{\bf aut}_{\mF_x} \mO_d \simeq K^x :=
\left\{  
g \in \ud \ : \ \wa g (t) = t \ , \ t \in \mF_x
\right\} \ .
\]
\noindent Moreover, suppose $\mO_\mSUE \subseteq \mF$ and that $\mF$ is generated by the intertwiners spaces $(\rhors)$, $r,s \in \bN$. Then, $K^x \subseteq \sud$ and $\mF_x \simeq \mO_{K^x}$ for every $x \in X$.
\end{lem}

\begin{proof}
Since $\mO_\mUE \subseteq \mF$ we find $\mF' \cap \coe \subseteq \mO'_\mUE \cap \coe = C(X)$, thus $\mF'_x \cap \mO_d = \bC 1$. Let now $\alpha \in {\bf aut}_{\mF_x} \mO_d$, $\psi , \psi' \in (\iota , H) \subseteq \mO_d$, $t \in \mF_x$. Then $\psi^* \alpha (\psi') t = \psi^* \alpha (\rho_x (t) \psi) = t \psi^* \alpha (\psi)$, and we conclude $\psi^* \alpha (\psi') \in \bC 1$. Thus $\alpha$ acts on $(\iota , H)$ as a unitary operator: $\alpha (\psi) = g \psi$, $\psi \in (\iota , H)$, $g \in \ud$, and this proves that $\alpha = \wa g$. We conclude that every element of ${\bf aut}_{\mF_x} \mO_d$ is the image of an element of $K^x$ w.r.t. the map (\ref{24}).

Let $\rho_x \in {\bf end} \mF_x$ be the canonical endomorphism $\rho_x (t_x) := 1 \otimes t_x$, $t_x \in (\rhors)_x \subset \mF_x$. Since $\mO_\mSUE \subseteq \mF \subseteq \coe$, we have \sC algebra inclusions $\mO_{\sud} \subseteq \mF_x \subseteq \mO_d$; with the argument of Rem.\ref{rem_dualpoint}, we conclude that $\mF_x$ is the fixed-point algebra of $\mO_d$ w.r.t. the action of the stabilizer $K^x$ of $\mF_x$ in $\mO_d$. Rem.\ref{rem_dualpoint} also implies that $K^x$ can be identified as a closed subgroup of $\sud$, acting on $\mO_d$ by the action (\ref{24}), so that there is an isomorphism $\mF_x \simeq \mO_{K^x}$.
\end{proof}

\begin{cor}
\label{lem_def_G_x}
Let $\mE \ra X$ be a rank $d$ vector bundle, $\mG \subseteq \mUE$ a closed group. Then, for every $x \in X$ the action (\ref{24}) induces an isomorphism
\[
{\bf aut}_{(\cog)_x} \simeq \mG^x := 
\left\{ 
g \in \ud \ : \ \wa g (t) = t \ , \ t \in (\cog)_x  
\right\}
\]
\noindent Moreover, if $\mG \subseteq \mSUE$ then $(\cog)_x = \mO_{\mG^x}$ for every $x \in X$.
\end{cor}

\begin{rem} {\em
Let $\mG \subseteq \mSUE$. By the previous corollary, $\cog$ is a continuous bundle with fibres nuclear (\cite{DLRZ02}), simple (\cite[Thm.3.1]{DR87}) \sC algebras.
} \end{rem}

The groups $\mG^x$ are said the {\bf spectral fibres} of $\mG$.

Now, recall that for every $u \in \mcUE$ there is $g \in \mUE$ such that $u = g_x$, $x := \eta (u)$ (see Rem.\ref{rem_full}). Thus, for every $x \in X$, $\mUE_x \simeq \ud$ acts by automorphisms on $(\coe)_x \simeq \mO_d$, by defining $\wa u (t) := \wa g_x (t)$, $t \in (\coe)_x$.

\begin{defn}
\label{def_group_bundle}
Let $\mE \ra X$ be a rank $d$ vector bundle, $\mG \subseteq \mUE$ a closed group. The {\bf spectral bundle} associated with $\mG$ is the topological subspace of $\mcUE$
\[
\mcG :=
\left\{  
u \in \mcUE \ : \ \wa u (y) = y  \ , \ y \in (\cog)_x \ , \ x \in X
\right\} 
\subseteq \mcUE \ ,
\]
\noindent endowed with the projection $\pi_\mG : \mcG \ra X $, $\pi_\mG (u) := \eta (u)$.
\end{defn}

By recalling the definition of spectral fibre, we obtain an isomorphism $\pi_\mG^{-1} (x) \simeq \mG^x$ for every $x \in X$. We denote by 
\begin{equation}
\label{def_sg}
\mSG :=  \left\{ g \in \mUE \ : \ g_x \in \mcG \  \forall x \in X \right\}
\end{equation}
\noindent the closed group of elements of $\mUE$ which are continuous sections of $\mcG$. Since $\mG_x \subseteq \mG^x$, $x \in X$, it is clear that 
\[
\mG \subseteq \mSG \ \ .
\]
\noindent Note that there is a simple criterion to determine whether some $g \in \mUE$ is an element of $\mSG$: it suffices in fact to verify that $g_x \in \mG^x$ for every $x \in X$.

\begin{ex} {\em
\label{ex_tsb}
We consider a closed group ${G_0} \subseteq \ud$, acting on the trivial bundle $\mE := X \times H$ by the natural action $v \mapsto g v$, $g \in {G_0}$, $v \in \mE_x \simeq H$. Note that $(\ers) = C(X) \otimes (\hrs)$, $r,s \in \bN$, and $\coe = C(X) \otimes \mO_d$. Let $\mG \subset \mUE = C(X,\ud)$ be the group of $G_0$-valued, constant maps, and $C(X,{G_0})$ the group of continuous maps from $X$ into $G_0$. It is clear that $C(X,{G_0})$ is a closed group of $\mUE$, and it is easily verified that
\[
(\ers)_{\mG} = (\ers)_{C(X,{G_0})} = C(X) \otimes (\hrs)_{G_0} \ .
\]
\noindent Thus $\mO_G \simeq C(X) \otimes \mO_{G_0}$; the spectral bundle associated with $\mG$ is $\mcG = X \times {G_0}$, with $\mSG = C(X,{G_0})$.
} \end{ex}

\bigskip

In analogy with \cite[Cor.3.3]{DR87}, the following result allows the reconstruction of the spectral bundle of $\mG$ from the triple $(\cog,\sigma_\mG,\mE)$ and thus, because of the amenability, from the dual ${\wa \mG}$. We denote by ${\bf aut}_{\cog} \coe$ the group of automorphisms of $\coe$ coinciding with the identity on $\cog$, i.e. the stabilizer of $\cog$ in $\coe$.

\begin{prop}
\label{dual_group_bundle}
For every closed group $\mG \subseteq \mUE$ there is a natural isomorphism $\mSG \simeq {\bf aut}_{\cog} \coe$, where $\mSG$ is the group of continuous sections of the spectral bundle $\mcG \ra X$. The {\bf duality map} $\left\{ \mcG \mapsto \wa \mG \right\}$ is injective.
\end{prop}

\begin{proof}
By definition it is clear that there is an immersion $\mSG \hra {\bf aut}_{\cog} \coe$, $g \mapsto \wa g$. Viceversa, let $\alpha \in {\bf aut}_{\cog} \coe$. Then for $\psi , \psi' \in \wE$ we find $\psi^* \alpha (\psi) y = \psi^* \alpha (\psi y) = \psi^* \sigma (y) \alpha (\psi) = y \psi^* \alpha (\psi)$ for every $y \in \cog$, so that by Prop.\ref{rel_com_og} $\psi^* \alpha (\psi)$ belongs to $C(X)$ and $\alpha = \wa g$ for some $g \in \mUE$. In order to prove that $g \in \mSG$ it suffices to verify, for every $x \in X$, that $g_x$ belongs to the spectral fibre $\mG^x$: $\wa g (y) = y$ $\Rightarrow$ $(\wa g (y))_x = \wa g_x (y_x) = y_x$, $y \in \cog$. Since $(\cog)_x = \mO_{\mG^x}$, we find $g_x \in \mG^x$. Finally, let $\mG$, $\mG'$ such that $\wa \mG = \wa {\mG'}$. Then $\cog = \mO_{\mG'}$, and $\mG^x = \mG'^x \simeq {\bf aut}_{(\cog)_x} \mO_d$, for every $x \in X$. Thus $\mcG = \mcG'$.
\end{proof}

The fact we recover $\mSG$ instead of $\mG$ in Prop.\ref{dual_group_bundle} is an example of the absence of the Galois property remarked in \cite[\S 7]{BL97}. A crucial point to get the injectivity of the duality map $\left\{ \mcG \mapsto \wa \mG \right\}$ is that the vector bundle $\mE$ has to be fixed; counterexamples will be given in the sequel (Cor.\ref{no_uni}).

\begin{cor}
\label{cor_iso_g}
Let $\mG , \mG' \subseteq \mUE$ with $\mG = \mSG$, $\mG' = SG'$. Suppose there exists an automorphism $\alpha$ of $\coe$ with $\alpha (\sigma^r,\sigma^s) = (\sigma^r,\sigma^s)$ for $r,s \in \bN$, and such that $\alpha (\cog) = \mO_{\mG'}$. Then $\mG$ and $\mG'$ are conjugate, i.e. there exists $u \in \mUE$ such that $\mG' = u \mG u^*$.
\end{cor}

\begin{proof}
Since $\alpha (\sigma^r,\sigma^s) = (\sigma^r,\sigma^s)$, by amenability we find $\alpha (\iota , \mE) = (\iota , \mE)$, thus $\alpha = \wa u$, $u \in \mUE$. Let now $\wa g \in {\bf aut}_{\cog} \coe \simeq \mG$, $\psi \in (\iota,\mE)_{\mG'}$. Since $\alpha^{-1} (\psi) = u^* \psi \in (\iota,\mE)_\mG$, we find $u g u^* \psi = \psi$; by the same argument, $(u g u^*) \wa{} \ $ is the identity on $(\ers)_{\mG'}$. Thus, $(u g u^*) \wa{} \in {\bf aut}_{\mO_{\mG'}} \coe \simeq \mG'$. By exchanging the role of $\mG$, $\mG'$ we obtain the desired equality $\mG' = u \mG u^*$.
\end{proof}

\begin{lem}
\label{lem_gsg}
Let $\mG \subseteq \mUE$ be a closed group. Then $\wa \mG = \wa{\mG'}$ for every closed group $\mG' \subseteq \mUE$ such that $\mcG = \mathcal {G}'$. In particular, $\wa \mG = \wa { \mSG }$  and $\cog = \mO_{\mSG}$.
\end{lem}

\begin{proof}
If $y \in (\ers)_\mG$, then $y_x \in (\cog)_x$ for every $x \in X$, thus by Cor.\ref{lem_def_G_x} $y_x$ is $\mG^x$-invariant. Let now $g \in \mSG$; since $g_x \in \mG^x$ for every $x \in X$ we find $(\wa g (y))_x = \wa g_x (y_x) = y_x$, and $\wa g (y) = y$. Now, if $\mG'$ has the same spectral bundle as $\mG$ then ${g'}_x \in \mG^x$ for every $g' \in \mG'$, thus $\wa {g'} (y) = y$. So that, we proved that $(\ers)_\mG \subseteq (\ers)_{\mG'}$. By exchanging the role of $\mG$, $\mG'$ we obtain $(\ers)_{\mG} = (\ers)_{\mG'}$.
\end{proof}

\begin{ex} {\em
\label{ex_ud}
Referring to Example \ref{ex_point}, it is clear that $\mG \neq C(X,\ud)$; anyway, the spectral bundle of $\mG$ is $X \times \ud$, and $\wa \mG = { C(X,\ud) }{\wa { }}$.
} \end{ex}

\bigskip

We now study some topological properties of spectral bundles, and investigate how the structure of the algebras $\cog$ reflects such properties. According to (\ref{def_lcud}), if $\pi_U : \mE |_U \ra U \times H$ is a local chart for $\mE$, we use the notation $\eta_U : \mcUE |_U \ra U \times \ud$ for the corresponding local chart induced on the unitary bundle.

\begin{defn}
Let $\mG \subseteq \mUE$ be a closed group. A {\bf local chart for the spectral bundle} $\mcG$ is given by a local chart $\pi_U : \mE |_U \ra U \times H$ such that $\eta_U (\mcG |_U) = U \times G_0$, where $G_0 \subseteq \ud$ is a closed group.
\end{defn}

Given a local chart $\pi_U$ for $\mcG$, it is clear that the restriction map $\pi_U^{1,1} : \mUE \ra C(U , \ud)$ restricts to a map $\pi_U^{1,1} : \mSG \ra C(U , G_0)$.

\bigskip

We say that the spectral bundle $\mcG \ra X$ is {\em locally trivial} if for every $x \in X$ there is a neighborhood $U \ni x$ defining a local chart for $\mcG$. If we can pick $U = X$, then we say that $\mcG$ is {\em trivial}. A subgroup $\mG$ of $\mUE$ does not define locally trivial spectral bundles in general (see Ex.\ref{ex_ord2} below). {\em In order for a more concise terminology, in the sequel we will say that $\mG$ is locally trivial if and only if the spectral bundle $\mcG \ra X$ is locally trivial}. If $X$ is connected, the local triviality implies that the isomorphism class of the spectral fibres is constant, so that there is $G_0 \subseteq \ud$ and an open trivializing cover $\left\{ U_i \right\}$ such that $\left. \mcG \right|_{U_i} \simeq U_i \times G_0$. Thus, the elements of $\mSG$ can be described in terms of continuous maps $g_i : U_i \ra G_0$ satisfying the cocycle relations $u_{ij} g_j = g_i u_{ij}$, where $u_{ij} : U_i \cap U_j \ra \ud$ are transition maps for $\mE$. If $\mcG$ is trivial, then there is an isomorphism $\mSG \simeq C (X , G_0)$.

Locally trivial group bundles of the type above have been studied in the setting of a generalized equivariant $K$-theory in \cite{NT04}.

\bigskip

\begin{ex} {\em
Let $\mE \ra X$ be a vector bundle. Then the spectral bundle associated with $\mUE$ is the unitary bundle $\mcUE$. In the same way, the spectral bundle associated with $\mSUE$ is the special unitary bundle $\mcSUE$. These bundles are locally trivial, and are trivial if and only if $\mE$ is the tensor product of a trivial bundle with a line bundle.
} \end{ex}

\begin{ex} {\em
\label{ex_ord2}
Let $X := [0,2]$, $\mE := X \times H$, $G_0 \subset \sud$ be a closed group. We fix $\omega \in (0,1)$ and consider the group
\[  
\mG = 
\left\{ 
 g \in C( X , \sud) : g_x \in G_0 \ , \ x \in [ 0 , \omega ] \
\right\} \subset \mSUE \ .  
\]
\noindent Since $\mE$ is trivial, we find $\coe \simeq C(X) \otimes \mO_d$. If $t \in \cog \subset \coe$, we find that $t_x \in \mO_d$ is $\sud$-invariant for every $x \in (\omega,2]$, i.e. $t_x \in \mO_\sud$. The map $\left\{ X \ni x \mapsto t_x \in \mO_d  \right\}$ being continuous, we obtain that $t_\omega$ is norm limit of elements of $\mO_\sud$, thus $t_\omega \in \mO_\sud$. This fact implies that $(\cog)_\omega = \mO_\sud$, so that $\mG^\omega = \sud$. Thus,
\[
\mcG = ([0,\omega) \times G_0) \  \sqcup\ ([\omega,2] \times \sud)
\]
\noindent is not locally trivial. This also implies that $(\ers)_\mG$, $r,s \in \bN$, is not the bimodule of continuous sections of a vector bundle: in fact, the Banach bundle associated with $(\ers)_\mG$ has fibre $(\hrs)_{G_0}$ for $x < \omega$, while for $x \geq \omega$ the fibre is $(\hrs)_\sud$. Note that $\mG_\omega = G_0$ is strictly contained in $\mG^\omega = \sud$. If $g \in \mSG$ then $g_x \in G_0$ for every $x \in [0,\omega)$; thus, $G_0$ being closed, by continuity we obtain $g_\omega \in G_0$. We conclude that $g \in \mG$, so that $\mG = \mSG$. 
} \end{ex}

\begin{ex} {\em
\label{ex1_ord2}
Let $S^1$ denote the circle, and $X := S^1 \times S^1$. Note that there is a standard closed cover $X = C_1 \cup C_2$, where $C_1$, $C_2$ are homeomorphic to the cylinder $[0,1] \times S^1$, and $Y := C_1 \cap C_2$ is the disjoint union $S^1 \sqcup S^1$. Let $\mE \ra X$ be a rank $d$ vector bundle, $d > 1$. Since $C_i$, $i = 1,2$, is homotopic to $S^1$, we find that $\mE |_{C_i}$ is trivial, so that $\mE$ is described by a transition map $u : Y \ra \ud$. An element $g \in \mSUE$ is described by a pair of continuous maps $g_i : C_i \ra \sud$, $i = 1,2$, satisfying the cocycle relation $g_2 = {\mathrm {ad}}u (g_1) := u g_1 u^*$. Let $G_0 \subseteq \sud$ be a closed subgroup, $\mG := \left\{ g \in \mSUE : g_2 (x) \in G_0 , x \in C_2 \right\}$; then the spectral bundle $\mcG$ is obtained by clutching the trivial bundles $C_1 \times \sud$, $C_2 \times G_0$ {\em via} the transition map ${\mathrm {ad}}u : Y \ra \left. \ud \right/ \bT$. It is clear that $\mcG$ is not locally trivial.
} \end{ex}

\begin{lem}
\label{lem_spbt}
If the spectral bundle $\mcG \ra X$ is trivial, then there is a compact Lie group $G_0 \subseteq \ud$ acting on $\mE$, in a such way that $(\ers)_\mG = (\ers)_{G_0}$, $r,s \in \bN$. If $\mG$ is locally trivial with spectral fibre $G_0$, then for every $x \in X$ there is a neighborhood $U \ni x$ such that $(\ers)_\mG |_U   = (\mE^r |_U , \mE^s |_U)_{G_0}$, $r,s \in \bN$.
\end{lem}

\begin{proof}
Let $G_0$ be the spectral fibre of $\mcG$. Then $\mcG \simeq X \times G_0$, and there is an isomorphism $\mSG \simeq C(X , G_0)$. Thus $\mcG$ can be associated with the group of constant $G_0$-valued sections (see Ex.\ref{ex_tsb}). The second assertion follows trivially from the first one, by considering local charts $\mcG |_U \simeq U \times G_0$, $U \subseteq X$.
\end{proof}

\begin{defn}
Let $\mG \subseteq \mUE$ be a closed group, $U \subseteq X$ a closed set with nonempty interior. A {\bf local chart for the pair} $(\coe , \cog)$ is given by 
\begin{itemize}

\item  a local chart $\alpha_U : \coe |_U \ra C(U) \otimes \mO_d$, such that 
\[
\alpha_U (\ers) |_U  =  C(U) \otimes ( H^r , H^s ) \  , \ r,s \in \bN \ ;
\]
\item  a closed group $G_0 \subseteq \ud$, such that $\alpha_U (\cog |_U) = C(U) \otimes \mO_{G_0}$.
\end{itemize}
\end{defn}

Note that if $\alpha_U$ is a local chart for $(\coe , \cog)$, then $\alpha_U (\ers)_\mG |_U  =  C(U) \otimes ( H^r , H^s )_{G_0}$, $r,s \in \bN$. Thus $\alpha_U$ 'trivializes' the Banach bimodules $(\ers)_\mG$, i.e. the spaces of arrows of $\wa \mG$. The pair $(\coe , \cog)$ is said {\em locally trivial} if for every $x \in X$ there is a closed $U \ni x$ with nonempty interior carrying a local chart for $(\coe , \cog)$.

\begin{prop}
\label{prop_lt}
Let $X$ be connected, $\mG \subseteq \mUE$ a closed group. Then, the pair $(\coe , \cog)$ is locally trivial if and only if $\mG$ is locally trivial. In such a case, there is a closed group $G_0 \subseteq \ud$ such that, for every $x \in X$,
\begin{equation}
\label{eq_fibre}
(\cog)_x \simeq \mO_{G_0}  \ \ , \ \  \mG^x \simeq G_0  \ \ .
\end{equation}
\end{prop}

\begin{proof}
Let $\pi_U : \mE |_U \ra U \times H$ be a local chart for $\mcG$, with the induced restriction morphism $\pi_U^{1,1} : \mSG \ra C(U , G_0)$. We now construct a local chart for $(\coe , \cog)$. Let
\[
\alpha_U : \coe |_U \ra C(U) \otimes \mO_d \ \ , 
\]
\noindent be the local chart defined by (\ref{def_flc}); we now verify that $\alpha_U$ restricts to a local chart for $\cog$. Let $g_0 \in G_0$; with an abuse of notation, we denote by $g_0$ the constant map $\left\{ U \ni x \mapsto g_0 \right\}$. It is clear that $g_0$ is a continuous (constant) section of $U \times G_0$. Since there is a bundle isomorphism $U \times G_0 \simeq \mcG |_U$, there is a continuous section $g_U : U \ra \mcG |_U$ such that 
\[
\pi_U^{1,1} (g_U) = g_0 \ \ .
\]
\noindent Since $(g_U)_x \in \mG^x$, $x \in U$, we find that $\wa g_U (y_U) = y_U$ for every $y_U \in \cog |_U$. Moreover, since $\alpha_U (g_U) = \pi_U^{1,1} (g_U) = g_0$, it is clear that 
\begin{equation}
\label{eq_af}
\alpha_U \circ \wa g_U  = \wa g_0 \circ \alpha_U \ \ .
\end{equation}
\noindent The previous equality implies that $\alpha_U (\cog |_U)$ is the fixed-point algebra of $C(U) \otimes \mO_d$ w.r.t. the action by automorphisms of the type $\wa g_0$, $g_0 \in G_0$, thus $\alpha_U (\cog |_U) = C(U) \otimes \mO_{G_0}$ (see Ex.\ref{ex_tsb}).

Let now $\alpha_U$ be a local chart for $(\coe , \cog)$; we denote by $\alpha_x \in {\bf aut} \mO_d$, $x \in X$ the automorphism induced by the evaluation of $\alpha_U$ over each fibre $x \in U$, so that there is a commutative diagram
\[
\xymatrix{
           \coe |_U
		    \ar[r]^-{\alpha_U}
		    \ar[d]
		 &  C(U) \otimes \mO_d
		    \ar[d]
		 \\ \mO_d
		    \ar[r]^-{\alpha_x}
		 &  \mO_d
}
\]
\noindent where the vertical arrows are the evaluation epimorphisms over the fibres. Since $\alpha_U (\iota , \mE) |_U = C(U) \otimes H$, the Serre-Swan equivalence implies that there is a local chart $\pi_U : \mE |_U \ra U \times H$. Now, the local chart $\pi_U$ induces a local chart $\eta_U : \mcUE |_U \ra U \times \ud$; in order to prove the proposition, we verify that $\eta_U$ restricts to a map from $\mcG$ onto $U \times G_0$. Let $g \in \mUE$ with $g_x \in \mG^x \subseteq \mcUE$, $(x , g'_x) := \eta_U (g_x) \in U \times \ud$. We verify that $g'_x \in G_0$. For this purpose, note that if $t \in (\ers)$ then, with the same argument as (\ref{eq_af}) and by recalling (\ref{24}),
\[
(\alpha_U \circ \wa g (t))_x = 
\alpha_x \circ \wa g_x (t_x) = 
\wa g'_x \circ \alpha_x (t_x) \ \ .
\]
\noindent From the previous equalities, it follows that $t_x \in (\cog)_x$ (i.e., $\alpha_x (t_x) \in \mO_{G_0}$) if and only if $\wa g'_x \in G_0$. Thus, the proposition is proved.
\end{proof}

\subsection{The cases $\mG = \mUE$, $\mG = \mSUE$.}

\begin{prop}
\label{lem21}
$\mO_{\mUE}$ is the $\sigma$-stable \sC subalgebra of $\coe$ generated by $\theta$ and $C(X)$, and is isomorphic to $C(X) \otimes \mO_{\ud}$.
\end{prop}

\begin{proof}
Let $\theta_d \in \mO_\ud$ denote the exchange operator on $H \otimes H$, $\sigma_d$ the canonical endomorphism on $\mO_{\ud}$. By \cite[Lemma 3.6]{DR87} it follows that $C(X) \otimes \mO_{\ud}$ is generated as a \sC algebra by $C(X)$, $1 \otimes \theta_d$, and by closing w.r.t. the action of $\iota \otimes \sigma_d$ (here $\iota$ denotes the identity automorphism on $C(X)$). Let $\left\{ U_i \right\}$ be a closed trivializing cover for $\mE$, with local charts $\pi_{U_i} : \mE |_{U_i} \ra {U_i} \times H$; we denote by $\left\{ u_{ij} := \pi_i \circ \pi_j^{-1} : X_i \cap X_i \ra \ud  \right\}$ the associated set of transition maps. By (\ref{def_flc}), local charts
\begin{equation}
\label{eq_lc}
\alpha_{U_i} : \coe |_{U_i} \ra C(U_i) \otimes \mO_d
\end{equation}
\noindent are induced, in such a way that the cocycle associated with $\coe$ as a \sC algebra bundle is given by $\left\{ \wa u_{ij} : X_i \cap X_j \ra {\bf aut} \mO_d \right\}$ (see \cite[\S 4]{Vas}). Since a local chart for $\mcUE$ is induced by every $\pi_{U_i}$, it follows from Prop.\ref{prop_lt} that (\ref{eq_lc}) restricts to a local chart $\beta_{U_i} : \mO_{\mUE} |_{U_i} \ra C(U_i) \otimes \mO_\ud$. Thus, the lemma will be proved if every transition map
\[
\beta_{ij} := \beta_{U_i} \circ \beta_{U_j}^{-1} 
\in {\bf aut} ( C(U_i \cap U_j) \otimes \mO_\ud )
\]
\noindent reduces to the identity. Now, it is clear that $\beta_{ij}$ is the restriction of $\wa u_{ij}$ to $C(U_i \cap U_j) \otimes \mO_\ud$; since $\theta_d$ is $\ud$-invariant, we find $\beta_{ij} ( 1 \otimes \theta_d ) = \wa u_{ij} ( 1 \otimes \theta_d ) = 1 \otimes \theta_d$; furthermore, $\beta_{ij} \circ (\iota \otimes \sigma_d) = \wa u_{ij} \circ (\iota \otimes \sigma_d) = (\iota \otimes \sigma_d) \circ \wa u_{ij} = (\iota \otimes \sigma_d) \circ \beta_{ij}$. Thus, by evaluating $\beta_{ij}$ over products of elements of the type $(\iota \otimes \sigma^r_d) (1 \otimes \theta_d)$, $r \in \bN$, we conclude that $\beta_{ij}$ is the identity automorphism, and the lemma is proved. In particular, for every $r \in \bN$ there is an isomorphism $(\mE^r , \mE^r)_{\mUE} \simeq ( H^r , H^r)_{\ud} \otimes C(X)$, while for $r \neq s$ we find $(\mE^r , \mE^s)_{\mUE} = \left\{ 0 \right\}$ (see \cite[Lemma 3.6]{DR87}).
\end{proof}

\begin{cor}
Let $\mE$, $\mE' \ra X$ be rank $d$ vector bundles. Then 
\[  
(\ers)_\mUE \simeq \delta_{r,s} \cdot \left[ ( \hrs)_{\ud} \otimes C(X) \right] \ ,
\]
\noindent where $\delta_{r,s}$ is the Kroneker symbol. There is an isomorphism of tensor \sC categories $\wa {\mUE} \simeq \wa { {\bf U} \mathcal {E'}}$ if and only if $\mE$, $\mE'$ have the same rank.
\end{cor}

\noindent The previous lemma implies that the \sC algebra $\mO_{\mUE}$ does not maintain geometrical informations about $\mE$, except for the rank. By restricting the group we can recover some geometrical data, as for example the first Chern class in the case of $\mSUE$. We are going to prove this fact. Let us consider the following endomorphism on $\mO_{\mUE}$, defined by the shift
\begin{equation}
\label{def_tau}
\rho (y) := P \otimes y \ ,
\end{equation}
\noindent where $P$ is the totally antisymmetric projection (\ref{23} , \ref{27}) and $y \in (\mE^r,\mE^r)_{\mUE}$. Note that the relation $\rho (y) = P \cdot \sigma_\mUE^d (y)$ holds for $y \in \mO_\mUE$.

\bigskip

Let $k \in \bN$; we denote by $\lambda \mE^k$ the tensor power of $\lambda \mE$ iterated $k$ times. For $k < 0$ we define $\lambda \mE^k := [(\lambda \mE)^* ]^{-k}$, where $(\lambda \mE)^*$ is the dual bundle of $\lambda \mE$. Note that $(\iota , \lambda \mE^k) = (\lambda \mE^{-k} , \iota)$ for every $k \in \bZ$. We denote by $\otimes_X$ the internal tensor product of Hilbert $C(X)$-bimodules.

\begin{prop}
\label{lem22}
$\mO_{\mSUE}$ is the $\sigma$-stable \sC subalgebra of $\coe$ generated by $\theta$ and $(\iota,\lambda \mE)$, and is isomorphic to the crossed product $\mO_{\mUE} \rtimes^{\lambda \mE}_\rho \bN$.
\end{prop}

\begin{proof}
By \cite[Lemma 3.7]{DR87} we find 
\[
(\hrs)_\sud = \delta_{r-s , dk} \cdot (H^r , H^r)_\ud \cdot R^k
\]
\noindent where $\delta_{i,j}$ is the Kroneker symbol, $R \in (\iota , H^d)$ is the totally antisymmetric isometry (\ref{21}), and $k$ is a positive integer (since $(\ers) = (\mE^s , \mE^r)^*$, it suffices to consider the case $k \geq 0$). Thus, by local triviality we conclude that $(\ers)_{\mSUE} \neq \left\{ 0  \right\}$ if and only if $r-s = kd$, $k \in \bZ$. In particular, $(\mE^r , \mE^r)_{\mSUE} = (\mE^r , \mE^r)_{\mUE}$, $r \in \bN$. Let now $\left\{ R_i  \right\}$ be the set of generators of $\lambda \mE$; then $\sum_i R_i^* R_i = 1$, so that $\sum_{I} R_{I}^* R_{I} = 1$, where $R_{I} := R_{i_1} \cdots R_{i_k} \in (\iota , \mE^{dk})$, $k \in \bN$. Let $y \in (\mE^r , \mE^{r+kd})_{\mSUE}$, $k > 0$; then $y = \sum_{I} ( y R_{I}^*) \ R_{I}$, with $( y R_{I}^*) \in (\mE^{r} , \mE^{r})_{\mSUE} = (\mE^{r} , \mE^{r})_{\mUE}$. The previous argument implies that 
\[
(\mE^r , \mE^{r+kd})_{\mSUE} = (\mE^r , \mE^r)_{\mSUE} \cdot (\iota , \lambda \mE^k) 
\ , \
k \in \bZ \ .
\]
\noindent Thus, we have proved that $\mO_{\mSUE}$ is generated as a \sC algebra by $\mO_{\mUE}$ and $(\iota , \lambda \mE)$. Furthermore, we have $R y = P \otimes y \cdot R \otimes 1 = \rho (y) R$, $R \in (\iota , \lambda \mE)$, $y \in (\ers)_{\mSUE}$, so that $\rho$ is inner in $\mO_{\mSUE}$, induced by $(\iota , \lambda \mE)$. In order to prove the lemma, we have to verify the universal property w.r.t. covariant representations $(\pi : \mO_\mUE \ra L(\mM_\pi) , \wa {\lambda \mE}_\pi )$. Let $\left\{ \psi \mapsto \psi_\pi  \right\}$, $\psi \in (\iota , \lambda \mE)$, $\psi_\pi \in \wa {\lambda \mE}_\pi$ denote the isomorphism between $(\iota , \lambda \mE) \subset \mO_\mSUE$ and $\wa {\lambda \mE}_\pi \subset L(\mM_\pi)$. Since 
\[
(\lambda \mE , \lambda \mE ) = P \cdot (\mE^d , \mE^d)_{\mUE} \cdot P = C(X) \cdot P \ ,
\]
\noindent for every $\psi , \psi' \in \lambda \mE$ there is $f_{\psi , \psi'} \in C(X)$ such that $\psi' \psi^* = f_{\psi , \psi'} P$. Furthermore, $\left \langle \psi , \psi' \right \rangle = \psi^* \psi' = \psi^*_\pi \psi'_\pi \in C(X)$. Now, $\psi' \psi^* \in \mO_{\mUE}$, so that if $\varphi \in (\iota , \lambda \mE)$ then
\begin{equation}
\label{eq_pp}
\begin{array}{ll}
\pi (\psi' \psi^*) \varphi_\pi & = f_{\psi , \psi'} \pi (P) \varphi_\pi \\ & =  
f_{\psi , \psi'} \varphi_\pi = (f_{\psi , \psi'} \varphi)_\pi = \\ & =
\psi'_\pi (\psi^* \varphi) = \psi'_\pi (\psi^*_\pi \varphi_\pi)  \ .
\end{array}
\end{equation}
\noindent Let $\left\{ R_{i,\pi}  \right\} \subset \wa {\lambda \mE}_\pi$, $\sum_i R_{i,\pi} R^*_{i,\pi} = \pi (P)$, be a set of generators for $\wa {\lambda \mE}_\pi$.  Then (\ref{eq_pp}) implies
\[
\pi (\psi' \psi^*) =  \pi (\psi' \psi^*) \pi (P) = 
\sum_i \pi (\psi' \psi^*)  R_{i,\pi} R^*_{i,\pi} = 
\sum_i \psi'_\pi \psi^*_\pi  R_{i,\pi} R^*_{i,\pi} \ ,
\]
\noindent so that $\pi (\psi' \psi^*) = \psi'_\pi \psi^*_\pi$. Thus, we define the $C(X)$-morphism 
\[
\Pi : \mO_{\mSUE} \ra L(\mM_\pi)  \ , \ \Pi ( t \varphi) := \pi (t) \varphi_\pi \ \ ,
\]
\noindent $t \in \mO_{\mUE}$, $\varphi \in (\iota , \lambda \mE)$. The previous remarks imply that $\Pi$ extends $\pi$ as desired.
\end{proof}

\begin{cor}
\label{no_uni}
Let $\mE$, $\mE' \ra X$ be rank $d$ vector bundles. Then 
\[  (\ers)_{\mSUE} \simeq \delta_{r-s , dk} \cdot \left[ (\mE^r , \mE^r)_\mUE \otimes_{X} (\iota , \lambda \mE^k) \right]  \  , \]
\noindent where $k \in \bZ$ and $\delta_{i,j}$ is the Kroneker symbol. So that there is an isomorphism of tensor \sC categories $\wa {\mSUE} \simeq \wa { {\bf SU} \mathcal {E'} }$ if and only if $\mE$, $\mE'$ have the same rank and first Chern class.
\end{cor}

\begin{ex} {\em
\label{ex_sue}
Let $X$ be a compact Hausdorff space with $H^2 (X , \bZ) = \bZ$, $\mE := \mL \oplus \mL^* \ra X$ a vector bundle, where $\mL^*$ denotes the dual of a nontrivial line bundle $\mL \ra X$. By construction $(\mE,\mE)$ is nontrival as a continuous bundle of matrix algebras, in fact $\mE$ is not the tensor product of a trivial bundle by a line bundle. So that, $\mSUE \neq C(X , \mathbb{SU}(2))$, in fact $(\mE,\mE)$ is generated as a $C(X)$-algebra by $\mSUE$. Despite that, by the previous corollary there is an isomorphism of tensor \sC categories $\wa {\mSUE} \simeq { C(X , \mathbb{SU}(2) ) } \wa{} \ $: in fact, $\lambda \mE = \mL \otimes \mL^* \simeq X \times \bC$, so that the first Chern class of $\mE$ vanishes and $(\iota , \lambda \mE) \simeq C(X) \cdot R$, where $R$ is defined by (\ref{21}). Thus, $\mO_\mSUE \simeq C(X) \otimes \mO_{\mathbb {SU} (2)}$.
} \end{ex}

The previous example shows that the injectivity of the duality map (Prop.\ref{dual_group_bundle}) is verified only if the embedding into the category of vector bundles is fixed. In fact, here $\wa \mSUE$ is exhibited as a dual of non-isomorphic groups, i.e. $\mSUE \simeq {\bf aut}_{\mO_\mSUE} \coe$ and $C(X , \mathbb{SU}(2) ) \simeq {\bf aut}_{\mO_\mSUE} (C(X) \otimes \mO_2)$. The fact we recover different groups depends on the choice of embedding $\wa \mSUE$ into the category of tensor powers of $\mE := \mL \oplus \mL^*$ (corresponding to the inclusion $\mO_\mSUE \hra \coe$), or $X \times \bC^2$ (corresponding to $\mO_\mSUE \hra C(X) \otimes \mO_2$).

\bigskip

Anyway, in general the spectral fibres $\left\{ \mG^x \right\}_{x \in X}$ of a closed group $\mG \subseteq \mUE$ do not depend on the choice of the vector bundle $\mE'$ realizing the embedding $\cog \hra \mO_{\mE'}$. Thus, {\em different embeddings of $\wa \mG$ into the category of vector bundles correspond to different topologies over $\bigcup_x \mG^x$, and define not necessarily isomorphic spectral bundles}.

\section{Noncommutative Pullbacks.}
\label{bim_pul}

\subsection{Some general properties of Hilbert bimodules.}
Let $\mA$ be a \sC algebra, $\mM$ a Hilbert $\mA$-bimodule. We consider the Banach $\mA$-bimodules $(\mrs)$, $r,s \in \bN$, introduced in Sec.\ref{preli}. If $\mA$ is unital then $(\mA , \mA) \simeq \mA$, and there is a natural isomorphism $(\mA , \mM) \simeq \mM$. The category $\mM^\otimes$ with objects the tensor powers of $\mM$ and arrows $(\mrs)$ is a semitensor \sC category in the sense of \cite[\S 2]{DPZ97}; roughly speaking, a tensor product is defined on the objects, while on the arrows just the operation of tensoring on the right by the identity arrow is admitted. This structure reflects the well-known fact that in general it is not possible to define in a consistent way the tensor product of right $\mA$-module operators (\cite[\S 13.5]{Bla}). Anyway, it is also possible to associate with $\mM$ a {\em tensor} \sC category, having objects the tensor powers of $\mM$, and arrows the sets of $\mA$-bimodule operators {\em commuting with the left $\mA$-action} :
\[
\mcB (\mrs) := \left\{ 
t \in (\mrs) \ : \  at = ta \ , \ a \in \mA
\right\} 
\]
\noindent (see \cite[\S 2]{DPZ97} for details: the operation of tensor product makes sense for elements of $\mB (\mrs)$). Let $\com$ denote the \sC algebra associated with $\mM$ (recall that if $\mA$ is unital and $\mM$ is finitely generated, then $\com$ is the CP-algebra of $\mM$). Following \cite[\S 3]{DPZ97}, we construct the \sC subalgebra $\bm$ of $\com$ generated by the Banach bimodules $\mcB (\mrs)$. Note that the relations $\mcB (\mrs) = (\mrs) \cap \mA'$ hold in $\com$; in particular, $\mB (\mA , \mA) \simeq M(\mA) \cap \mA' = ZM(\mA)$ (the last equality is a classical remark due to Busby). Thus, $\bm \subseteq \mA' \cap \com$. The \sC algebra $\bm$ is naturally endowed with the canonical endomorphism 
\begin{equation}
\label{def_ctau}
\tau \in {\bf end} \bm  \ : \ \tau (t) := 1 \otimes t \ 
\in \mcB ( \mM^{r+1} , \mM^{s+1} ) \ \ , 
\end{equation}
\noindent where $t \in \mcB (\mrs)$ and $1$ is the identity on $\mM$. Note that in general $\tau$ cannot be extended to an endomorphism of $\com$, and $\tau$ is not the identity on elements of $\mB (\mA , \mA)$.

The next lemma allows to compute the relative commutant of $\mA$ and $\bm$ in relevant particular cases. We recall that by construction $\com$ ($\bm$) is endowed with a $\bZ$-grading, arising from the canonical circle action. We say that a \sC algebra $\mB \subseteq \com$ is $\bZ$-{\em graded} if $\cup_{k \in \bZ} (\mB \cap \com^k)$ is dense in $\mB$.

\begin{lem}
\label{normality}
Let $\mA$ be a unital \sC algebra, $\mM$ a Hilbert $\mA$-bimodule, $\mB \subseteq \bm$ a $\tau$-stable, $\bZ$-graded \sC algebra. Suppose there is an isometry $R \in \mB(\mA,\mM^d) \cap \mB$, $d \in \bN$, $d > 1$. Then, $\bm = \mA' \cap \com$. Furthermore, if there is $\lambda \in (-1,0) \cup (0,1)$ such that the equality $R^* \tau (R) = \lambda 1$ holds, then $\mA = \mB' \cap \com$ (and $\mA$ is normal in $\com$ in the sense of \cite[\S 1]{DPZ97}).
\end{lem}

\begin{proof}
The first assertion is proven in \cite[Prop.3.4]{DPZ97}. The second assertion follows from \cite[Prop.3.5]{DPZ97}, by defining the sequence $R_k := \tau^{k-1} (R) \cdots \tau (R) R$.
\end{proof}

\begin{rem} {\em
\label{rem_sym}
Let $\mA := C(X)$ be abelian, $\mM$ a finitely generated, projective Hilbert $C(X)$-bimodule. Then, $C(X)$ is a unital \sC sunbalgebra of the centre of $\bm$; moreover, by projectivity there is $n \in \bN$ with inclusions of Banach $C(X)$-bimodules $\mB (\mrs) \hra C(X) \otimes \bM_{n^r , n^s}$, $r,s \in \bN$. This implies that $\bm$ is a continuous bundle of \sC algebras over $X$. Let $\mE \ra X$ be the vector bundle having module of continuous sections isomorphic to $\mM$ as a right Hilbert $C(X)$-module ($\mE$ is the so-called {\bf symmetrization} of $\mM$ as in \cite[Def.1.5]{AE01}); we denote by $\sigma$ the canonical endomorphism of $\coe$. Now, there is an inclusion of \sC algebra bundles $i : \bm \hra \coe$; note that in general it is false that $\sigma \circ i = i \circ \tau$: the obstruction is given by the eventuality that $\tau (f) \neq f$ for some $f \in C(X)$.
} \end{rem}

We now pass to consider group actions over a Hilbert bimodule $\mM$, and the associated \sC algebra. Let
\[  
\mUM := \left\{ u \in \mcB ( \mM , \mM ) : uu^* = u^*u = 1 \right\} \ ;
\]
\noindent then, $\mUM$ acts by automorphisms on $\com$, by extending the map
\begin{equation}
\label{alpha_action}
\psi \mapsto \wa u (\psi) := u \psi \quad , \quad \psi \in \mM \  , \ u \in \mUM
\end{equation}
\noindent  or, equivalentely, by the analogue of (\ref{24}) for elements of $(\mrs)$. Note that $\bm$ is stable w.r.t. the action (\ref{alpha_action}). By the same argument used for the action (\ref{24}), the canonical endomorphism $\tau$ commutes with $\wa u$, $u \in \mUM$. Since $\mUM \subset \mcB (\mM , \mM)$, we find $\wa u (a) = uau^* = a$, $u \in \mUM$, $a \in \mA$; thus $\mUM$ acts on $\com$, $\bm$ by $\mA$-bimodule automorphisms.

\begin{defn}
\label{defgm}
Let $\mM$ be a Hilbert $\mA$-bimodule, $\mG \subseteq \mUM$ a closed group. We denote by ${\wa \mG}_\mM$ the semitensor \sC category with objects the tensor powers of $\mM$ and arrows the invariant $\mA$-bimodules
\begin{equation}
\label{mrsg}
(\mrs)_\mG := \left\{ 
t \in (\mrs) : t = \wa g (t) := g^{\otimes^s} \cdot t \cdot g^{* \ \otimes^r} , g \in \mG 
\right\}  \ .
\end{equation}
\end{defn}

We denote by $\cmg$ the \sC subalgebra of $\com$ generated by the invariant spaces $(\mrs)_\mG$. The following definition encodes a particular class of group actions over Hilbert bimodules.

\begin{defn}
\label{dual_bimod_act}
Let $\mA$ be a \sC algebra, $\mM$ a Hilbert $\mA$-bimodule. A {\bf tensor} $\mG${\bf -action} over $\mM$ is given by a closed group $\mG \subseteq \mUM$ such that $(\mrs)_\mG \subseteq \mcB (\mrs)$, $r,s \in \bN$.
\end{defn}

Let $\mG \subseteq \mUM$ be a closed group acting on $\mM$ by a tensor action. Then, the following elementary properties hold:

\begin{enumerate}

\item  the condition $\mA \subseteq (\mA , \mA)_\mG \subseteq \mB ( \mA , \mA) = ZM (\mA)$ forces $\mA$ to be an abelian \sC algebra;

\item  since $(\mrs)_\mG \subseteq \mcB (\mrs)$, $r,s \in \bN$, we obtain that ${\wa \mG}_\mM$ is a {\em tensor} \sC category, and $\cmg \subseteq \bm$;

\item  $\cmg$ is $\tau$-stable (in fact, $\tau \circ \wa u = \wa u \circ \tau$, $u \in \mUM$). We denote by $\tau_\mG \in {\bf end} \cmg$ the corresponding restriction.

\end{enumerate}

\subsection{Basic properties of noncommutative pullbacks.}
As remarked in \cite[Rem.6.4]{BL97}, in general a well-defined symmetry in the sense of Doplicher and Roberts fails to exist in $\mM^\otimes$: the exchange operators $\theta (r,s)$ introduced in \S \ref{cavb_gb} are well defined as elements of $(\mM^{r+s},\mM^{r+s})$ if and only if the left and right $\mA$-module actions coincide, i.e. $a \varphi = \varphi a$, $a \in \mA$, $\varphi \in \mM$. In fact
\[
\theta (r,s) (( \psi \otimes \psi') a ) = 
\theta (r,s) ( \psi \otimes (\psi' a)) = 
(\psi' a) \otimes \psi = 
\psi' \otimes (a \psi) \ ,
\]
\noindent while
\[
(\theta (r,s) ( \psi \otimes \psi' )) a = 
(\psi' \otimes \psi)  a = 
\psi' \otimes (\psi a) \ ,
\]
\noindent $\psi \in \mM^r$, $\psi' \in \mM^s$. We now introduce a class of Hilbert bimodules such that the corresponding semitensor \sC categories of tensor powers admit a 'maximal' symmetric tensor \sC subcategory. We start with two preliminary remarks.

\begin{rem} {\em
Let $\mA$ be a \sC algebra, $\mM$ a Hilbert $\mA$-bimodule. We introduce the following unital, abelian \sC algebra:
\begin{equation}
\label{def_xm}
C(X_\mM) := \left\{ f \in ZM(\mA) : f a \psi a' =  a \psi fa' \ , \psi \in \mM , a,a' \in \mA   \right\} \ .
\end{equation}
\noindent Since the identity of $M(\mA)$ belongs to $C(X_\mM)$, it is clear that $\mA$ is a $C(X_\mM)$-algebra. In the case in which $\mA$ is unital, we find $f \in C(X_\mM)$ $\Leftrightarrow$ $f \psi =  \psi f$, $\psi \in \mM$.
} \end{rem}

\begin{rem} {\em
Let $\mA$ be a $C_0(X)$-algebra, $\mN$ a Hilbert $C_0(X)$-bimodule. Then, the algebraic tensor product $\mN \odot_{C_0(X)} \mA$ with coefficients in $C_0(X)$ is endowed with a natural $\mA$-valued scalar product $\left \langle \psi \otimes a , \psi' \otimes a' \right \rangle := \left \langle \psi , \psi' \right \rangle \cdot a^* a'$, $\psi , \psi' \in \mN$, $a,a' \in \mA$. We denote by $\mN \otimes_X \mA$ the corresponding completition. $\mN \otimes_X \mA$ is a right Hilbert $\mA$-module in the natural way.
} \end{rem}

\begin{defn}
\label{def_bp}
Let $\mA$ be a $C_0(X)$-algebra, $\mM$ a Hilbert $\mA$-bimodule such that $C_0(X) \subseteq C(X_\mM)$. $\mM$ is called a {\bf noncommutative pullback} ({\bf nc-pullback}, in the sequel) if there is a vector bundle $\mE \ra X$ with an isomorphism of {\bf right} Hilbert $\mA$-modules $\mM \simeq \wE \otimes_X \mA$. $\mM$ is said {\bf full} if $X$ is compact and $C(X) = C(X_\mM)$.
\end{defn}

\noindent Every nc-pullback is generated as a right Hilbert $\mA$-module by elements of $\wE$. By using the formalism of amplimorphisms (in the sense of \cite[\S 1]{DR89}), nc-pullbacks correspond in the unital case to \sC algebra morphisms of the type $\phi : \mA \ra \mA \otimes \bM_d$ such that $\phi(1) \in C(X) \otimes \bM_d$ and $\phi (f) = \phi (1) f$, $f \in C(X)$. The corresponding module is recovered as $\mM := \left\{ \psi \in \bC^d \otimes \mA : \phi(1) \psi = \psi \right\}$, with right $\mA$-action given by the scalar multiplication and left $\mA$-action $a , \psi \mapsto \phi (a) \psi$.

\begin{rem} {\em
\label{rem_full_ncp}
Let $\mA$ be a $C(X)$-algebra (with $C(X)$ unital), $\mM \simeq \wE \otimes_X \mA$ a nc-pullback. Then, $C(X)$ is a unital \sC subalgebra of $C(X_\mM)$, so that there is a surjective map $p : X_\mM \ra X$. We consider the pullback bundle $\mE_\mM := \mE \times_X X_\mM \ra X_\mM$ (in the sense of \cite[I.1.16]{Kar}), so that $\wE_\mM \simeq \wE \otimes_X C(X_\mM)$. Now, a natural isomorphism $(\wE \otimes_X C(X_\mM)) \otimes_{X_\mM} \mA \simeq  \wE \otimes_X \mA$ of right Hilbert $\mA$-modules is induced by the map $(\psi \otimes f) \otimes a \mapsto \psi \otimes (fa)$, $\psi \in \wE$, $f \in C(X_\mM)$, $a \in \mA$. Thus, in the case in which $C(X)$ is unital and $\mA$ is a $C(X)$-algebra, every nc-pullback is a full nc-pullback.
} \end{rem}

\begin{rem} {\em
When $X$ is compact and $\mA$ is unital and separable, a nc-pullback $\mM$ naturally defines a class $(\mM , 0)$ in the $KK$-theory group $\mR KK( X ; \mA,\mA)$ (see \cite{Kas88} about the previous notation); thus $\mM$ is a central $\mA$-bimodule in the sense of \cite{PT00}. Let us now denote by $\mR KK_f ( X ; \mA,\mA)$ the group of Kasparov bimodules which are finitely generated as right Hilbert $\mA$-modules. Then, there is a forgetful morphism $\pi : \mR KK_f ( X ; \mA,\mA) \ra K_0(\mA)$, defined by assigning to each Kasparov bimodule the corresponding right Hilbert $\mA$-module. If $\mA = C(Y)$ is commutative there is a surjective map $p : Y \ra X$, and every vector bundle $\mE' \ra Y$ defines a Kasparov bimodule $(\wa {\mE'} , 0) \in \mR KK( X ; C(Y) , C(Y) )$. Thus $\pi : \mR KK_f ( X ; C(Y) , C(Y) ) \ra K^0 (Y)$ is an epimorphism. If $\mE \ra X$ is a vector bundle, every nc-pullback of the type $\mM \simeq \wE \otimes_X C(Y)$ defines a class $(\mM , 0)$ in $\pi^{-1} [p_* \mE]$, where $p_* \mE \ra Y$ is the pullback bundle. Thus $\pi^{-1} [p_* \mE] \subseteq \mR KK_f ( X ; C(Y) , C(Y) )$ classifies the nc-pullbacks of $\mE$ over $C(Y)$.
} \end{rem}

\begin{ex} {\em
Let $\mM$ be a Hilbert $\mA$-bimodule isomorphic as a right Hilbert $\mA$-module to the free module $\bC^d \otimes \mA$. Then, $\mM$ is a nc-pullback of the trivial rank $d$ vector bundle over $X_\mM$. Bimodules of this type appear in the framework of so-called 'Hilbert \sC systems' as in \cite{BL97,BL01}, and are called {\em algebraic Hilbert spaces}. Other examples, given by Hilbert bimodules arising from contractions in compact metric spaces, can be found in \cite[\S 4]{PWY01}.
} \end{ex}

\begin{ex} {\em
Let $X$ be a compact Hausdorff space, $\mZ$ a unital, abelian $C(X)$-algebra, $\mE \ra X$ a vector bundle, $\rho \in {\bf end}_X \mZ$ a $C(X)$-endomorphism. Then, we can define a nc-pullback $\mM := \wE \otimes_X \mZ$, with left $\mZ$-action $z \psi := \psi  \rho (z)$, $\psi \in \mM$, $z \in \mZ$. The case in which $\rho$ is the identity corresponds to the usual notion of pullback of a vector bundle.
} \end{ex}

\begin{ex} {\em
Let $\mA$ be a \sC algebra with identity $1_\mA$, $E \in C(X) \otimes \bM_d$ a projection defining a vector bundle $\mE \ra X$. Let $\mB$ be the corner $(E \otimes 1_\mA) \cdot (C(X) \otimes \bM_d \otimes \mA) \cdot (E \otimes 1_\mA)$. Then, every unital $C(X)$-morphism $\phi : C(X) \otimes \mA \ra \mB$ corresponds with a nc-pullback of $\mE$, isomorphic to the external tensor product $\wE \otimes \mA$ as a right Hilbert $(C(X) \otimes \mA)$-module, and defining a natural structure of vector $\mA$-bundle in the sense of \cite{MF80}.
} \end{ex}

\bigskip

Let $X$ be a compact Hausdorff space, $\mA$ a unital $C(X)$-algebra with identity $1$, $\mM \simeq \wE \otimes_X \mA$ a nc-pullback. According to Rem.\ref{rem_full_ncp}, we assume that $\mM$ is full, so that we identify $X$ with $X_\mM$. The following elementary properties hold:

\begin{enumerate}

\item There is a map $j : \wE \ra \mM$, $j(\psi) := \psi \otimes 1$. It is clear that $j(\wE)$ is contained in $\mM \equiv (\mA , \mM)$. Thus, $j(\wE)$ is a finitely generated Hilbert $C(X)$-bimodule in $\com$ with support $1$, and is isomorphic to $\wE$. 

\item By universality of the CP-algebra, $j$ extends to a $C(X)$-monomorphism $j : \coe \hra \com$.

\item The relation $j (\ers) \subset (\mrs)$ holds, for every $r,s \in \bN$ (in particular, $j (\wE) = j (\iota , \mE) \subset (\mA , \mM)$). So that, $j$ defines an injective functor $j : \mE^\otimes \hra \mM^\otimes$ of semitensor \sC categories. It is clear that $j$ restricts to an isomorphism $j : \mE^\otimes \ra j(\wE)^\otimes$.

\item By amenability we have $\wa \sigma \simeq \mE^\otimes$, where $\sigma \in {\bf end}\coe$ is the canonical endomorphism (Prop.\ref{rel_com_og}); so that, there is an injective functor $j' : \wa \sigma \hra \mM^\otimes$ of semitensor \sC categories.

\end{enumerate}

\begin{lem}
With the above notation, suppose that $\bm' \cap \com = \mA$. Then:

\begin{enumerate}

\item $\com \cap \com' = C(X) \simeq C(X_\mM)$;

\item an inner endomorphism $\sigma_\mE \in {\bf end}\com$ is induced by $j(\wE)$;

\item $\sigma_\mE (t) = \tau (t)$, $t \in \bm$;

\item $(\sigma_\mE^r , \sigma_\mE^s ) = j( \ers ) \subset (\mrs)$, $r,s \in \bN$;

\item there is an isomorphism $\wa \sigma \simeq \wa \sigma_\mE$ of tensor \sC categories, and $\sigma_\mE$ has permutation symmetry.

\item For every $r,s \in \bN$, the equality
\[
       (\mrs) =  \sigma_\mE^s (\mA) \cdot j( \ers ) = 
       j( \ers ) \cdot \sigma_\mE^r (\mA) 
\]
\noindent holds; i.e., $\left\{ \sigma_\mE^s (z) t = t \sigma_\mE^r (a)  ,  z \in \mZ , t \in j( \ers ) \right\}$ is total in $(\mrs)$.

\end{enumerate}
\end{lem}

\begin{proof}
\
\begin{enumerate}

\item  The elements of $C(X)$ commute with $\mA$ and $\mM \equiv (\mA , \mM)$ in $\com$, thus $C(X) \subseteq \com \cap \com'$. Viceversa, $\com' \cap \com \subseteq \bm' \cap \com  =  \mA$, thus $\com' \cap \com \subseteq \mA \cap \mA'$; so that, if $f \in \com' \cap \com$ then $f \psi = \psi f$, $\psi \in \mM$, and $f \in C(X)$.

\item  $j(\wE)$ is a finitely generated Hilbert $C(X)$-bimodule in $\com$ with support $1$; thus, the inner endomorphism $\sigma_\mE \in {\bf end} \com$ induced by $j(\wE)$ is defined, and the point is proved. Note that $j(\wE) = ( \iota , \sigma_\mE)$, $j(\wE)^r := {\mathrm{span}} \left\{ \right. \psi_1 \cdots \psi_r , \psi_i \in j(\wE) , i = 1, \ldots ,r  \left. \right\} = (\iota , \sigma_\mE^r)$, $r \in \bN$.

\item  Let $t \in \mB (\mrs)$, $\psi \in j(\wE) \subseteq (\mA , \mM)$. Then, $(\psi \otimes 1_s) \cdot  t = \psi \otimes t = (1 \otimes t ) \cdot ( \psi \otimes 1_r )$, where $1_r$ is the identity of $(\mM^r \ \mM^r)$. Thus, the relation $\psi t = \tau (t) \psi$ holds; on the other hand, by definition $\psi t = \sigma_\mE (t) \psi$.

\item  By definition, $j (\ers) \subseteq (\mrs)$ and $j(\iota , \mE^r) = j(\wE)^r$. Now, $j( \ers )$ is generated as a Banach $C(X)$-bimodule by elements of the type $\psi_\vL \psi_\vM^*$, $\psi_\vL := \psi_{l_1} \cdots \psi_{l_s} \in j(\wE)^s$, $\psi_\vM \in j(\wE)^r$ (see \S \ref{preli}). Since $j(\wE)^r = (\iota , \sigma_\mE^r)$, we find $j( \ers ) \subseteq (\sigma_\mE^r , \sigma_\mE^s )$. Viceversa, if $t \in (\sigma_\mE^r , \sigma_\mE^s )$ then $t = \sum_{LM} \psi_\vL t_{LM} \psi_\vM^*$, where $t_{LM} := \psi_\vL^* t \psi_\vM \in$ $(\sigma^s , \iota) \cdot ( \sigma_\mE^r , \sigma_\mE^s ) \cdot (\iota , \sigma_\mE^r)$ $\subseteq (\ii) = C(X)$ ($\iota \in {\bf end} \com$ denotes the identity automorphism). Thus, $t \in j( \ers )$.

\item  By the previous point, we find $j(\wE)^\otimes \simeq \wa \sigma_\mE$. Moreover, by amenability we obtain $\mE^\otimes \simeq \wa \sigma$. Since $\mE^\otimes \simeq j(\wE)^\otimes$, we conclude $\wa \sigma \simeq \wa \sigma_\mE$. Since $\sigma$ has permutation symmetry (Rem.\ref{rem_ps}), we obtain that $\sigma_\mE$ has permutation symmetry.

\item  Since $\sigma_\mE^r (a) = \sum_M \psi_\vM a \psi_\vM^* \in (\mM^r , \mM^r)$, $a \in \mA$, the inclusion $j( \ers ) \cdot \sigma_\mE^r (\mA) \subseteq (\mrs)$ is proved. Viceversa, if $t \in (\mrs)$ then $t = \sum_{LM} \psi_\vL t_{LM} \psi_\vM^*$, where $\psi_\vL \in j(\wE)^s$, $\psi_\vM \in j(\wE)^r$, $t_{LM} := \psi_\vL^* t \psi_\vM \in$ $(\mM^s , \mA) \cdot (\mrs) \cdot (\mA , \mM^r) \subseteq \mA$. Note that $\psi_\vL t_{LM} = \sigma_\mE^s (t_{LM}) \psi_\vL$, $t_{LM} \psi_\vM^* = \psi_\vM^* \sigma_\mE^r (t_{LM})$.

\end{enumerate}

\end{proof}

\bigskip

The next lemma (generalizing the second statement of Lemma \ref{normality}) will be used in the sequel, and allows to compute the relative commutant of certain $\tau$-stable $\bZ$-graded subalgebras of $\bm$.

\begin{lem}
\label{normality2}
With the above notation, let $\mB \subseteq \bm$ be a $\tau$-stable, $\bZ$-graded \sC algebra. Suppose there is a finitely generated Hilbert $C(X)$-bimodule $\mR \subseteq j(\iota,\mE^d) \cap \mB$, $d > 1$, with support $P_\mR$, such that 
\[
\mR \mR^* := {\mathrm{span}} \left\{ R' R^* \ , \ R,R' \in \mR  \right\} =  C(X) \cdot P_\mR  \  .
\]
\noindent Moreover, suppose that for some $\lambda \in (-1,0) \cup (0,1)$ the equality $R^* \tau (R') = \lambda R^* R'$ holds for every $R,R' \in \mR$. Then $\mA = \mB' \cap \com$.
\end{lem}

\begin{proof}
The equality $\mR \mR^* = C(X) \cdot P_\mR$ and the Serre-Swan Theorem imply that $\mR$ is the module of continuous sections of a line bundle $\mL \ra X$, $\mL \subset \mE^d$. Let us consider a closed cover $\left\{ X_i \right\}$ of $X$, trivializing $\mE , \mL$. Now, $C(X)$ is a unital \sC subalgebra of the centre of $\com$; thus $\com$ is a $C(X)$-algebra, that we regard as un upper semicontinuous bundle over $X$. 

Let $\mA_i$ denote the restriction of $\mA$ over $X_i$ as an upper semicontinuous bundle. The restriction $\com |_{X_i}$ of $\com$ over $X_i$ is generated as a \sC algebra by $\mM_i := \mM \otimes_\mA \mA_i$, so that $\com |_{X_i} \simeq \mO_{\mM_i}$; we denote by $\alpha_i : \com \ra \mO_{\mM_i}$ the associated epimorphism, so that $\alpha_i (\mA , \mM) = (\mA_i , \mM_i)$. Since $\mE |_{X_i} \simeq X_i \times \bC^d$ is trivial, we find that $\mM_i$ is isomorphic as a right Hilbert $\mA_i$-module to the free module $\bC^d \otimes \mA_i$. In particular, $\mR_i := \alpha_i (\mR)$ is generated as a Hilbert $C(X_i)$-bimodule in $\mO_{\mM_i}$ by an isometry $R_i \in (\mA_i , \mM_i^d)$, corresponding to the generator of the free module of continuous sections of ${\mL} |_{X_i} \simeq X_i \times \bC$. Note that $R_i^* \tau (R_i) = \lambda \alpha_i(1)$, where $\alpha_i (1)$ is the identity of $\mO_{\mM_i}$. By using $R_i$ as in Lemma \ref{normality}, we construct a sequence $\left\{ R_{i,k} \right\}$.

Let now $t \in \mB' \cap \com$; then, every $t_i := \alpha_i(t)$ commutes with elements of $\left\{ R_{i,k} \right\}$, so that $t_i \in \mA_i$. We now consider a partition of unity $\left\{ \lambda_i \right\} \subset C(X)$ subordinate to $\left\{ X_i \right\}$; then, $\lambda_i t_i = \lambda_i t \in \mA$, and we conclude that $t = \sum_i \lambda_i t_i \in \mA$.
\end{proof}

\bigskip

Let now $\mZ$ be an abelian, unital $C(X)$-algebra, $\mM \simeq \wE \otimes_X \mZ$ a full nc-pullback. We consider a closed group $\mG \subseteq \mSUE$ such that
\[
\mG = \mSG \ ,
\]
\noindent  where $\mSG$ is defined by (\ref{def_sg}). Recall that there is a monomorphism $j : \coe \hra \com$ such that $j(\ers) \subset (\mrs)$; in particular, $j (\mUE) \subset (\mM , \mM)$. Let us now suppose that $j(\mG)$ acts on $\mM$ by a tensor action (recall Def.\ref{dual_bimod_act} and subsequent remarks). In order for more concise notations, we identify $\mG$ with $j(\mG)$ and define
\[
\mB := \cmg \subseteq \bm \ \ , \ \ \rho := \tau_\mG \in {\bf end} \mB \ ;
\]
\noindent we also consider the group 
\[
{\bf aut}_\mB (\com , \sigma_\mE) :=
\left\{ 
\alpha \in {\bf aut} \com \ : \
\alpha |_{\mB} = \iota \ , \  
\alpha \circ \sigma_\mE = \sigma_\mE \circ \alpha
\right\} \ .
\]
\noindent  Our purpose is to study the \sC dynamical systems $(\mB , \rho)$, $(\com , \mG)$. In the particular case in which $\mZ = C(X)$ (i.e. $\mM = \wE$), the \sC dynamical system $(\mB , \rho)$ is of the type $(\cog , \sigma_\mG)$ studied in the previous section.

\begin{prop}
\label{cross_bp}
With the above notation, the following properties hold:

\begin{enumerate}
\item $\mB' \cap \mB = \mB' \cap \com = \mZ$ (thus, $\bm' \cap \com = \mZ$, $\com \cap \com' = C(X)$).
\item The map $\left\{ \mG \ni g \mapsto \wa g \right\}$ defines an isomorphism $\mG \simeq {\bf aut}_\mB (\com , \sigma_\mE)$.
\end{enumerate}
\end{prop}

\begin{proof}
\
\begin{enumerate}

\item There are obvious inclusions $\mZ \subseteq \mB' \cap \mB \subseteq \mB' \cap \com$. Moreover, we recall Lemma \ref{lem_scp} and apply Lemma \ref{normality2} with $\mR := (\iota,\lambda \mE) \subset (\iota , \mE^d)_\mG \subset \mB \subseteq \bm$, $\lambda := (-1)^{d-1} d^{-1}$. Thus $\mB' \cap \com = \mZ$, and the others equalities immediately follow.

\item Let $g \in \mG$. Since $\mG$ acts on $\mM$ by a tensor action, $\wa g \in {\bf aut} \com$ restricts to the identity on $\mB$; moreover, $(\iota , \sigma_\mE) = j(\wE)$ is $\mG$-stable, thus $\wa g \circ \sigma_\mE = \sigma_\mE \circ \wa g$ and $\wa g \in {\bf aut}_\mB (\com , \sigma_\mE)$. Viceversa, let $\alpha \in {\bf aut}_{\mB} (\com , \sigma_\mE)$. If $\psi , \psi' \in (\iota , \sigma_\mE) = j(\wE)$ we find $\psi^* \alpha (\psi') \in \mB' \cap \com$. By the previous point, $\psi^* \alpha (\psi') \in \mZ$. Let now $\theta \in (\mE^2 , \mE^2)_\mG$ be the operator defined by (\ref{defsim})). Since $\alpha \circ \sigma_\mE = \sigma_\mE \circ \alpha$, we find $\rho ( \psi^* \alpha (\psi') ) = \psi^* j(\theta) \cdot \sigma_\mE \circ \alpha (\psi') = \psi^* j(\theta) \alpha (j(\theta) \psi') = \psi^* \alpha (\psi')$. Thus $\psi^* \alpha (\psi)$ is $\rho$-invariant; since $\rho = \sigma_\mE |_\mB$, we conclude that $\psi^* \alpha (\psi)$ commutes with elements of $(\iota , \sigma_\mE)$. Since $\com$ is generated as a \sC algebra by $\mZ$, $(\iota , \sigma_\mE)$, we find $\psi^* \alpha (\psi') \in \com \cap \com' = C(X)$. Thus $(\iota , \sigma_\mE) \simeq \wE$ is $\alpha$-stable, and  $\alpha = \wa g$ for some $g \in \mUE$. Now, we can regard at $\wa g$ as an automorphism of $\coe$ leaving $\cog \subseteq \cmg$ pointwise invariant; thus, we apply Prop.\ref{dual_group_bundle} and conclude that $g \in \mG = \mSG$.

\end{enumerate}
\end{proof}

\begin{rem} {\em
Let $z \in \mZ$. It follows from the point 1 of the previous proposition that $\rho (z) =\sigma_\mE (z) = z$ $\Leftrightarrow$ $[z , (\iota , \sigma_\mE)] = 0$ $\Leftrightarrow$ $z \in C(X) = \com \cap \com'$.
} \end{rem}

\begin{cor}
\label{cross_bp1}
Let $\mG$ be locally trivial. Then, for every $r,s \in \bN$,
\begin{equation}
\label{gps_mrsg}
(\mrs)_\mG = \rho^s (\mZ) \cdot j(\ers)_\mG = j(\ers)_\mG \cdot \rho^r (\mZ)
\end{equation}
\noindent (i.e., the set $\left\{ \rho^s (z) t = t \rho^r(z) \ , \ z \in \mZ , t \in (\ers)_\mG \right\}$ is total in $(\mrs)_\mG$). Moreover, 
\begin{equation}
\label{gps_mrsg1}
(\mrs)_\mG = (\rhors) \ ,
\end{equation}
\noindent so that there is an isomorphism of tensor \sC categories $\wa \rho \simeq \wa \mG_\mM$.
\end{cor}

\begin{proof}
\
\begin{enumerate}

\item The inclusion $\rho^s (\mZ) \cdot j(\ers)_\mG \subseteq (\mrs)_\mG$ is trivial. Viceversa, let $\left\{ \psi_l \right\}$ be a finite set of generators of $j(\wE)$. The elements of $(\mrs)_\mG$ are of the type $t := \sum_{\vL \vM} \rho^s (t_{\vL \vM}) \psi_\vL \psi_\vM^*$, where $\psi_\vL := \psi_{l_1} \cdots \psi_{l_s} \in j(\iota , \mE^s)$, $\psi_\vM \in j(\iota,\mE^r)$, $t_{\vL \vM} := \psi_\vL^* t \psi_\vM \in {\mB}' \cap \com = \mZ$. Note that $\psi_\vL \psi_\vM^* \in j(\ers)$, thus $t \in \rho^s (\mZ) \cdot j(\ers)$. Let $G_0$ be the spectral fibre of $\mG$, and $U \subseteq X$ a local chart for $\mcG$. Then, $G_0$ acts on $\com |_U$ by automorphisms; if $t_U$ is the restriction of $t$ over $U$, by averaging w.r.t. the Haar measure of $G_0$ we obtain
\[  t_U = 
    \int_{G_0} \wa g_0 (t) \ dg_0  |_U = 
    \sum_{\vL \vM} \rho^s (t_{\vL \vM}) \int_{G_0} \wa g_0 (\psi_\vL \psi_\vM^* ) \ dg_0 
    \ |_U
    \ \ ;
\]
\noindent by Lemma \ref{lem_spbt}, we find
\[
\int_{G_0} \wa g_0 (\psi_\vL \psi_\vM^* ) \ dg_0 \ |_U \ \in \ 
( j(\ers)|_U )_{G_0} \ = \ j(\ers)_{\mG} |_U \ \ ,
\]
\noindent thus $t_U \in \ \rho^s (\mZ) \cdot (\ers)_{\mG} |_U$. By considering an open cover $\left\{ U_i \right\}$ trivializing $\mcG$, and a subordinate partition of unity $\left\{ \lambda_i \right\}$, we conclude that $t = \sum_i \lambda_i t_{U_i} \in \rho^s (\mZ) \cdot (\ers)_\mG $.

\item Let $t \in (\mrs)_\mG$. Then, by the previous point $t \in \rho^s(\mZ) (\sigma_\mE^r,\sigma_\mE^s) \cap \mB \subseteq (\rho^r,\rho^s)$. Viceversa, if $t \in (\rho^r,\rho^s)$, then $t_{\vL \vM}$ defined as above belongs to $\mZ$, thus $t \in (\mrs)$; in particular, since $t$ is $\mG$-invariant, we conclude $t \in (\mrs)_\mG$. The proof of the isomorphism $\wa \rho \simeq \wa \mG_\mM$ goes through the same line of Cor.\ref{cor_itc}.
\end{enumerate}
\end{proof}

\begin{rem} {\em
Let $\mE$ be a $G_0$-vector bundle, where $G_0$ is a locally compact group acting trivially on $X$. We denote by $\mG \subseteq \mUE$ the associated closed group in the sense of Rem.\ref{rem_glc}. Suppose that $\mM$ is a nc-pullback of $\mE$ carrying a tensor action by $\mG \subset \mSUE \subset \mUM$. Then, (\ref{gps_mrsg},\ref{gps_mrsg1}) hold. In fact, we can use the Haar measure of $G_0$ and apply the invariant mean argument used in the proof of (\ref{gps_mrsg}).
} \end{rem}

%
%

\bigskip
\bigskip

{\bf Acknowledgments.} The author would like to thank S. Doplicher, for having proposed (and successively supported with improvements and encouragements) the initial idea of this work for his PhD Thesis, and P. Zito, R. Conti, G. Morsella for precious help and stimulating discussions.


\end{document}